\documentclass[11pt]{article}
\usepackage{graphicx}
\usepackage{amsmath,latexsym,amsbsy,amssymb, color}
\usepackage{graphicx}
\def\argmin{\hbox{arg}\!\min}
\def\ds{\displaystyle}
\def\forall{\hbox{for all}~}
\def\L{{\bf L}}
\def\dint{\int\!\!\!\int}

\def\bfX{{\bf X}}

\def\ve{\varepsilon}

\def\X{{\cal X}}

\def\R{I\!\!R}

\def\implies{\Longrightarrow}
\def\vp{\varphi}

\def\vs{\vskip 2em}

\def\v{\vskip 1em}

\def\O{{\cal O}}

\def\M{{\cal M}}

\def\C{{\cal C}}

\def\Gd{{\cal G}_\delta}

\def\ov{\overline}

\def\Tilde{\widetilde}
\def\Hat{\widehat}
\def\bega{\begin{array}}
\def\enda{\end{array}}
\def\begi{\begin{itemize}}
\def\endi{\end{itemize}}
\def\meas{\hbox{meas}}

\def\bel{\begin{equation}\label}
\def\eeq{\end{equation}}
\def\sqr#1#2{\vbox{\hrule height .#2pt
\hbox{\vrule width .#2pt height #1pt \kern #1pt
\vrule width .#2pt}\hrule height .#2pt }}
\def\square{\sqr74}
\def\endproof{\hphantom{MM}\hfill\llap{$\square$}\goodbreak}

\newtheorem{theorem}{Theorem}[section]

\newtheorem{lemma}{Lemma}[section]
\newtheorem{proposition}{Proposition}[section]
\newtheorem{remark}{Remark}[section]
\newtheorem{definition}{Definition}[section]
\newtheorem{example}{Example}[section]
\newtheorem{conjecture}{Conjecture}[section]
\begin{document}

\title{\bf Generic Properties of First Order \\
Mean Field Games}
\vs

\author{Alberto Bressan$^{(*)}$ and 
Khai T.~Nguyen$^{(**)}$\\~~\\
 {\small $^{(*)}$~Department of Mathematics, Penn State University,}
 \\ {\small $^{(**)}$~Department of Mathematics, North Carolina State University.} \\~~\\
{\small E-mails: axb62@psu.edu,~khai@math.ncsu.edu.}
}
\maketitle

\begin{abstract} We consider a class of deterministic mean field games, 
where the state associated with each player evolves according to an ODE which is linear 
w.r.t.~the control.  Existence, uniqueness, and stability of solutions are 
studied from the point of view of generic theory.  Within a suitable topological space of 
dynamics and cost functionals, we prove that, for ``nearly all" mean field games
 (in the Baire category sense) the best reply map is single valued for a.e.~player.
As a consequence, the mean field game admits a strong (not randomized) solution.
Examples are given of open sets of games admitting a single solution, and other open sets
admitting multiple solutions.  Further examples show 
the existence of an open set of MFG
having a unique solution which is asymptotically stable w.r.t.~the best reply map, and
another open set of MFG having a unique solution which is unstable. We conclude with an example 
of a MFG with terminal constraints which does not have any solution, 
not even in the mild sense with randomized strategies.

\end{abstract}

\v

\section{Introduction}
\setcounter{equation}{0}
\label{s:1}
This paper deals with a class of mean field games  with a continuum of players, 
where the state associated with each player evolves according to a controlled ODE.
We study the existence, uniqueness, and stability of solutions from the point of view of generic theory.
Namely, we seek properties of solutions that are satisfied either on some open set of MFG, or
for ``nearly all" MFG in the topological sense \cite{Dug, W}; i.e., for all MFG in the intersection of 
countably many open dense sets.

Let $(\Omega,{\cal B}, \mu)$ be a probability space.  More precisely, we
assume that 
$\Omega$ is a metric space with Borel $\sigma$-algebra ${\cal B}$, while 
$\mu$ is an atomless probability measure on $\Omega$. Without loss of generality, throughout 
the following we assume $\Omega=[0,1]$ with Lebesgue measure.
We regard $\xi\in \Omega$ as a Lagrangian variable, labelling one particular player.
Accordingly,  we shall denote by $t\mapsto x(t,\xi)$ a  trajectory for player $\xi$.
By selecting one trajectory $x(\cdot ,\xi)\in \C\bigl([0,T];\,\R^n\bigr)$ 
for each player (depending measurably on $\xi$), one obtains an element $X$
in the space
\bel{LOC} \L^1\Big( \Omega\,;~\C\bigl([0,T];\, \R^n\bigr)\Big).\eeq
The space (\ref{LOC}) is naturally endowed with the Banach norm
\bel{NLOC}
\|X\|~\doteq~\int_\Omega \left( \sup_{t\in [0,T]} \bigl|x(t,\xi)\bigr|\right)\, d\xi.\eeq
To define a  (deterministic) mean field game, for each player $\xi\in \Omega$ we consider an 
optimal control problem where the dynamics and the cost functions
also depend on the cumulative distribution $X$ of all other players. 
To express this dependence, we consider a finite number of smooth 
scalar functions $\phi_1,\ldots,\phi_N
 \in \C^2\bigl([0,T]\times\R^n\bigr)$,
 and define
$\eta(t)=(\eta_1, \ldots,\eta_N)(t)$ to be the  vector of  ``moments"
\bel{etai} \eta_i(t)~=~\int_\Omega \phi_i \bigl(t, x(t,\xi)\bigr)  \,d\xi,\qquad i=1,\ldots,N.\eeq
The control problem for player $\xi$  takes the form
\bel{mmin}
\hbox{minimize:}\qquad \int_0^T L\bigl(t, x(t), u(t), \eta(t)\bigr)\, dt 
+ \psi\bigl(x(T)\bigr),
\eeq
subject to the dynamics 
\bel{mdyn} \dot x(t)~=~f\bigl( t, x(t), u(t), \eta(t)\bigr)\qquad\qquad t\in [0,T],\eeq
 and with initial datum
\bel{iitt}
x(\xi,0) \,=\,\bar x(\xi)
 .\eeq

\begin{definition}\label{d:21} In the above setting, by
a {\bf strong solution} to the mean field game  we mean a family 
of control functions $t\mapsto u(t,\xi)\in\R^m$ and corresponding trajectories $t\mapsto x(t,\xi)
\in\R^n$, 
defined for $\xi\in \Omega$ and $t\in [0,T]$,
such that the following holds.

For a.e.~$\xi\in \Omega$, the control $u(\cdot,\xi)$ and the trajectory
$x(\cdot,\xi)$ provide an optimal solution to the optimal control problem 
(\ref{mmin})--(\ref{iitt}) for player $\xi$,  
where $\eta(t)=(\eta_1,\ldots, \eta_N)(t)$ is the vector of moments defined at 
(\ref{etai}).
\end{definition}

A mean field game
thus yields a (possibly multivalued) map $\eta\mapsto \Phi(\eta) $ from $\C([0,T] ;\,\R^N)$ into itself.
Namely, given $\eta(\cdot)$, for each $\xi\in \Omega$ consider an optimal trajectory
$x^\eta(\cdot , \xi)$ of the corresponding optimal control problem (\ref{mmin})--(\ref{iitt}).
We then set 
\bel{Phi}
\Phi(\eta)~\doteq~\Tilde \eta ~=~(\Tilde\eta_1,\ldots, \Tilde \eta_N), 
\qquad\qquad  \Tilde\eta_i(t)~\doteq~\int_\Omega \phi_i\bigl(t, x^\eta(t,\xi)\bigr)\, d\xi,\eeq
under suitable assumptions that will ensure that the integral in (\ref{Phi}) is well defined.
By  definition, a fixed point of this composed map
\bel{mot}\bega{c}\qquad \eta(\cdot)\quad \mapsto\qquad \bigl\{ x^\eta(\cdot, \xi)\,;~\xi\in \Omega\bigr\}
\quad\mapsto\quad \Tilde\eta~=~\Phi(\eta)\\[2mm]
[\hbox{moments}]~\mapsto~ [\hbox{optimal trajectories}]~\mapsto~ [\hbox{moments}]
\enda
\eeq
yields a strong solution to the mean field game.

\begin{remark} {\rm
In general,  the map $\Phi$ can be multivalued.   
Indeed, for some $\eta(\cdot)$,
 there can be a subset   $V\subseteq \Omega$
with positive measure, such that each player  $\xi\in V$ has two or more optimal 
trajectories.
For this reason, a mean field game may not have a solution in the strong sense considered
in Definition~\ref{d:21}. In order to achieve a general existence theorem
one needs to relax the concept of solution, allowing the possibility of randomized
strategies \cite{BaF, CC, CP}.   This leads to the problem of finding a fixed point of an upper semicontinuous
convex-valued multifunction, which exists by Kakutani's theorem~\cite{Ce,K}.
}
\end{remark}

Following  the standard literature on fixed points of 
continuous or multivalued maps, we introduce
\begin{definition}\label{d:22}  A solution $x=x^*(t,\xi)$ to the above mean field game is {\bf stable}
if the corresponding function $\eta^*\in \C^0\bigl([0,T];\, \R^N\bigr)$ at (\ref{etai})
is a stable fixed point of the multifunction $\Phi$ at (\ref{Phi}).
Namely, for every $\ve>0$ there exists $\delta>0$ such that the following holds.
For every sequence 
$\bigl(\eta^{(k)}
\bigr)_{k\geq 0}$ such that 
\bel{etk}
 \|\eta^{(0)}-\eta^*\|_{\C^0}~<~\delta,\qquad \quad \eta^{(k)}~\in~\Phi\bigl(\eta^{(k-1)}\bigr)
 \quad\forall k\geq 1,
 \eeq
one has $ \|\eta^{(k)}-\eta^*\|_{\C^0}<\ve $ for all  $k\geq 1$.

If, in addition, every  such sequence $\bigl(\eta^{(k)}\bigr)$ converges to $\eta^*$, then we say that the solution is
{\bf asymptotically stable}.

If the solution is not stable, we say that it is {\bf unstable}.
\end{definition}

Next, we say that a solution of the mean field game is structurally stable if  it
persists under small perturbations of the dynamics and the cost functionals. More precisely:

\begin{definition}\label{d:23}  We say that a solution $x=x(t,\xi)$ 
to the above mean field game (\ref{etai}--(\ref{iitt}) is {\bf structurally stable} (or equivalently: {\bf  essential})
if, given $\ve>0$,  there exists $\delta>0$ such that the following holds.
For any 
perturbations $(f^\dagger, L^\dagger, \psi^\dagger, \phi^\dagger, \bar x^\dagger)$ satisfying
\bel{sper}\max\Big\{
\|f^\dagger-f\|_{\C^2}\,,~\|L^\dagger-L\|_{\C^2}\,,~
\|\psi^\dagger-\psi\|_{\C^2}\,,~
\|\phi^\dagger-\phi\|_{\C^2}\Big\} ~<~\delta,\qquad
\|\bar x^\dagger -\bar  x\|_{\L^\infty}<\delta,\eeq
the corresponding perturbed game
has a solution $x^\dagger =x^\dagger (t,\xi)$ such that 
\bel{stab1
}\sup_{t\in [0,T]} \int_\Omega \bigl| x^\dagger(t,\xi)- x(t,\xi) \bigr|\, d\xi~<~\ve.\eeq
\end{definition}

Throughout the following, we shall assume
that
the dynamics is affine w.r.t.~the control variable:
\bel{faff}f(x,u,\eta)~=~f_0(x,\eta) +\sum_{i=1}^m f_i(x,\eta) u_i\,,\eeq
and all functions $f,\psi,L$ have at least $\C^2$ regularity.

Since our MFG at (\ref{etai})--(\ref{iitt}) is characterized by the  5-tuple of functions $(f,L,\psi,\phi, \bar x)$,
we are interested in properties which are satisfied either
(i) for all games where 
 $(f,L,\psi,\phi, \bar x)$ ranges inside an open set (in a suitable Banach space), or
(ii) for {\bf generic} games, i.e., for all games where  $(f,L,\psi,\phi, \bar x)$ ranges over the 
intersection of countably many open dense sets.
Roughly speaking, 
the main results of the paper can 
be summarized as follows.

\begi 
\item[(i)] {\it Given a triple $(f,L,\phi)\in \C^3\times \C^3\times \C^3$, for a generic  pair $(\psi,\bar x)
\in \C^3\times \L^\infty$, the best reply map $\eta\mapsto\Phi(\eta)$ 
in (\ref{mot}) is 
single valued.    As a consequence, the MFG (\ref{etai})--(\ref{iitt}) admits a strong solution.}
\item[(ii)]  {\it There is an open set of mean field games with  a unique solution,
which is stable and essential.}
\item[(iii)]   {\it There is an open set of mean field games with  a unique solution,
which is unstable, and essential.}
\item[(iv)]  {\it  There is an open set of mean field games with two solutions, both essential.}
\endi
More precise statements of these results will be given in the following sections.
The remainder of the paper is organized as follows.

As a warm-up, in Section~\ref{s:3} we review the basic tools for proving generic properties. 
Here we consider a family of optimal control problems where the dynamics is linear w.r.t.~the control functions.
We show that, for generic dynamics $f$, running cost $L$ and terminal cost $\psi$, for a.e.~initial datum $x(0)=\ov x$  
the optimal control is unique.

Section~\ref{s:5} provides a simple way to construct mean field games with multiple solutions.
Given an optimal control problem and a pair $(x^*,u^*)$ 
(not necessarily optimal) which satisfies  the Pontryagin necessary conditions, we show the existence of
a  mean field game where $u^*$ is the optimal control for every  player.
As a consequence,  for any control problem
where the Pontryagin equations have multiple solutions, one can construct a MFG with multiple 
solutions. Under generic assumptions, all of these solutions are structurally stable.

Section~\ref{s:7} contains the main result of the paper. Namely, for a generic MFG
of the form (\ref{etai})--(\ref{iitt}), the best reply map $\eta\mapsto \Phi(\eta)$ is single valued.  
Hence the MFG admits a strong solution.  Here the analysis is far more delicate than in the proof of
the generic uniqueness for the optimal control problem in Section~\ref{s:3}.   Indeed, we need to show that
the statement
\begi
\item {\it 
The set of initial points $\ov x$, for which the problem (\ref{mmin})--(\ref{iitt}) has multiple  solutions,
has measure zero}
\endi
is true not just for one function $\eta(\cdot)$, but simultaneously  for all  
functions $\eta=(\eta_1,\ldots, \eta_N)$, in a suitable domain.

Finally, Section~\ref{s:2} collects a variety of examples, where the MFG have multiple strong solutions, 
Some of these are stable,  in the sense of Definition~\ref{d:22}, while others are unstable.

We conclude with two examples of MFG without solution. The first one is a well known case where 
nonexistence is due to the fact
that the best reply of each player is not unique. No strong solution exists, but one can construct 
a {\bf mild solution} where each player adopts a randomized strategy.
In the second example, the presence of a terminal constraint lacking a transversality condition
prevents the existence of any solution, even in the mild (randomized) sense.

Some concluding remarks, pointing to future research directions, are given in Section~\ref{s:9}.

Mean field games  with stochastic dynamics have been introduced by Lasry and Lions \cite{LL}  
and by Huang, Malham\'e and Caine \cite{HMC}, to model the behavior of a large number of interacting agents.
Their solution leads to a well known system of forward-backward parabolic equations.
Solutions to first order MFG (with deterministic dynamics) can be obtained as a vanishing viscosity limit
of these parabolic PDEs, i.e., as viscosity solutions to a corresponding Hamilton-Jacobi equation
\cite{CC, CCC, CG, CP}.
Equivalently, one can take a Lagrangian approach, describing the optimal control and the
optimal trajectory of each single agent.   This is the approach followed in the present paper.
Some examples of MFG with unique or  with multiple solutions can be found in~\cite{BaF}.
A concept of structural stability for solutions to first order MFG was proposed in \cite{BriCa}.

\section{Generic uniqueness for optimal control problems}
\setcounter{equation}{0}
\label{s:3}
%
Consider an optimal control problem of the form
\bel{oc1}
\hbox{minimize:}\qquad J[u]~\doteq~\ds\int_0^T L\bigl( x(t), u(t)\bigr)\, dt + \psi\bigl(x(T)\bigr),
\eeq
with dynamics which is affine  in the control:
\bel{oc2} \dot x(t)~=~f\bigl(x(t), u(t)\bigr)~=~f_0(x(t))+\sum_{i=1}^mf_i(x(t))\, u_i(t),\qquad \quad x(0)= \bar x.\eeq
Here  $u(t)\in \R^m$ while $x(t)\in \R^n$.  
 To fix ideas, we shall consider the couple $(f,L)$ satisfying the following assumptions.

\begi
\item[{\bf (A1)}] {\it 
 The functions $f_i:\R^n\mapsto \R^n$, $i=0,\ldots,m$, are  twice continuously differentiable. Moreover the vector fields $f_i$ satisfy the sublinear growth condition
 \bel{sublin}
 \bigl|f_i(x)\bigr|~\leq~c_1\, \bigl(|x|+1\bigr)
 \eeq
 for some constant $c_1>0$ and all $ x\in\R^n$.
\item [{\bf (A2)}] The running cost  $L:\R^n\times\R^m\mapsto\R$ is twice continuously differentiable and satisfies
 \bel{Lbig}\left\{
 \bega{rl}
 L(x,u)&\geq~c_2\bigl( |u|^2-1\bigr),\\[2mm]
 |L_x(x,u)|&\leq~\ell(|x|)\cdot (1+|u|^2),
 \enda\right.
\eeq
 for some constant $c_2>0$ and some continuous function $\ell$. 
 Moreover,  $L$ is uniformly convex w.r.t.~$u$. Namely, for some $\delta_L>0$, 
the $m\times m$ matrix of second derivatives w.r.t.~$u$ satisfies
\bel{u-convex-L}
L_{uu}(x,u)~>~\delta_L\cdot \mathbb{I}_m \qquad\qquad\forall x,u.
\eeq
Here $ \mathbb{I}_m$ denotes the $m\times m$ identity matrix
}
\endi

Throughout the following, the open ball centered at the origin with radius $r$ is denoted by
$B_r =B(0,r)$, while $\ov B_r$ denotes its closure. 
Under the previous assumptions, optimal controls and  optimal trajectories of  the optimization
problem (\ref{oc1})-(\ref{oc2}) satisfy uniform a priori bounds:

\begin{lemma}\label{bb} Assume that  the couple $(f,L)$ satisfies {\bf (A1)-(A2)} and $\psi:\R^n\to [0,\infty[$ is twice continuously differentiable. Then there exist continuous functions $\alpha,\beta:[0,\infty[\to [0,\infty[$ such that  the following holds. Given any  initial point $\bar{x}\in \ov B_r$, 
let  $u^*(\cdot)$ be an optimal control and let $x^*(\cdot)$  be the corresponding  optimal trajectory and  for the 
problem (\ref{oc1})-(\ref{oc2}).   Then 
\bel{prio}\hbox{\rm ess-}\!\!\!\sup_{t\in [0,T]}  |u^*(t)|~\leq~\alpha(r),\qquad\qquad 
\sup_{t\in [0,T]}|x^*(t)|~\leq~\beta(r).
\eeq
\end{lemma}
{\bf Proof.} Fix  $\bar{x}\in \ov B_r$. Calling $x_0(\cdot)$  the solution  of (\ref{oc2}) with $u(t)\equiv 0$, 
by (\ref{sublin}) it follows
\[
\sup_{t\in [0,T]}|x_0(t)|~\leq~(r+1)\cdot e^{c_1t}-1.
\]
Let $(x^*,u^*)$ be a pair of optimal trajectory and optimal control of the optimization
problem (\ref{oc1})-(\ref{oc2}). By the first inequality in (\ref{Lbig}), one has
\begin{multline}\label{u-L2}
\int_{0}^{T}\big|u^*(t)|^2dt~\leq~{1\over c_2}\,\left(\int_{0}^{T}L(x_0(t),0)dt+\psi(x_0(T))\right)+T\\
~=~{1\over c_2}\cdot\left(T\cdot \sup_{|y|\leq (r+1)\cdot e^{c_1T}-1} L(y,0)+\sup_{|y|\leq (r+1)\cdot e^{c_1T}-1}|\psi(y)|\right)+T~\doteq~\beta_1(r).
\end{multline}
Since   $x^*$ solves (\ref{oc2}) with $u\equiv u^*$, we have 
\[
|\dot{x}(t)|~\leq~c_1\cdot (|x|+1)\cdot \left(1+\sum_{i=1}^{m}|u^*_i(t)|\right)~\leq~{c_1\over 2}\cdot 
\bigl(|x|+1\bigr)\bigl(|u^*(t)|^2+m+2\bigr).
\]
Therefore, from (\ref{u-L2}) one obtains 
\[
\sup_{t\in [0,T]}|x^*(t)|~\leq~(r+1)\cdot \exp\left({c_1\over 2}\cdot [\beta_1(r)+(m+2)T]\right)-1~\doteq~\beta(r).
\]
To derive a pointwise bound on $u^*$, for every $\alpha\geq 0$ we consider the truncated function
\[
u_{\alpha}(s)~=~\begin{cases}
u^*(s)&\mathrm{if}\qquad |u^*(s)|\leq\alpha,\\
0&\mathrm{if}\qquad |u^*(s)|>\alpha.
\end{cases}
\]
Calling $x_{\alpha}$  the solution of (\ref{oc2}) with $u\equiv u_{\alpha}$, we have 
\[
\sup_{t\in [0,T]}|x_{\alpha}(t)|~\leq~\beta(r),\qquad \sup_{t\in [0,T]}|x^*(t)-x_{\alpha}(t)|~\leq~\beta_2(r)\cdot \int_{I_{\alpha}}|u^*(s)|ds
\]
for some continuous function $\beta_2$. For any constant $\gamma\geq 1$, setting 
$I_{\gamma}\doteq \bigl\{s\in [0,T]: |u^*(s)|>\gamma\bigr\}$ we estimate the difference in the costs:
$$\bega{l}\ds
0~\leq~J[u_{\gamma}]-J[u^*]~=~\int_{0}^{T}L(x_\gamma(t),u_\gamma(t))-L(x^*(t),u^*(t))dt+\psi(x_{\gamma}(T))-\psi(x_{*}(T))\\[4mm]
\qquad\ds
~\leq~\left((T+\beta_1(r))\cdot\sup_{|s|\leq \beta(r)}\ell(s)+\sup_{|y|\leq \beta(r)}|\nabla\psi(y)|\right)\cdot \beta_2(r)\cdot\int_{I_{\gamma}}|u^*(s)|ds\\[4mm]
\qquad\qquad \ds +\int_{I_{\gamma}}L(x^*(t),0)-L(x^*(t),u^*(t))dt\\[4mm]
\qquad\ds
~\leq~\alpha_1(r)\cdot  \int_{I_{\gamma}}|u^*(s)|ds-c_2\cdot \int_{I_{\alpha}}|u^*(s)|^2ds~\leq~\left(\alpha_1(r)-c_2\cdot\gamma\right) 
\enda$$
for some continuous function $\alpha_1(\cdot)$. This yields the first inequality in (\ref{prio}), with $\ds \alpha(r)={\alpha_1(r)/ c_2}$.
\endproof
In the following,  the positive cone in the Banach space $\C^2$ is denoted by
\bel{pos-cone}
\C^2_+(\R^n) ~\doteq~\left\{ \psi\in \C^2(\R^n)\,;~\inf_{x\in\R^n}\psi(x)> 0\right\}.
\eeq
We can now state the first result.
\begin{theorem} {\bf (Generic uniqueness for optimal control problems).} Under the assumptions {\bf (A1)-(A2)}, there  exists  a  $\Gd$ subset $\M\subset  \C^2_+(\R^n)$ such that the following holds.
For every $\psi\in\M$, the set of initial
points $\bar x\in \R^n$, for which   the optimal control problem (\ref{oc1})-(\ref{oc2})
has multiple solutions, has  Lebesgue measure zero.
\end{theorem}

{\bf Proof.}  
{\bf 1.} For every $\psi\in \C^2_+(\R^n)$  and  any integer $\nu\geq 1$, we consider a set 
of initial points yielding two distinct solutions:
\bel{Snu}
\bega{rl} S_\nu(\psi)&\doteq~\Big\{\bar x\in\ov B_{\nu},~
\hbox{the optimization
problem (\ref{oc1})-(\ref{oc2}) has two solutions}\\[2mm]
&  x_1(\cdot), ~x_2(\cdot), \ds\hbox{ with the same minimum cost,
and  with $\ds \bigl|x_1(T)-x_2(T)\bigr|\geq {1\over\nu}$ }\Big\}.
\enda
\eeq
Next, consider the set of terminal cost  $\psi$ leading  to a small set of multiple solutions:
\bel{key-1}
\M_\nu~\doteq~\left\{ \psi\in \C^2_+(\R^n)\,;~~\meas\bigl(S_\nu(\psi)\bigr)~<~{1\over\nu}\right\}.
\eeq
The theorem will be proved by showing  that $\M_{\nu}$ is open and  dense in $\C^2_+(\R^n)$. Indeed, if this is the case then  the set $\M =\ds\bigcap_{\nu\geq 1}\M_{\nu}$ is  a $\Gd$ subset of $\C^2_+(\R^n)$. Moreover, for any  $\psi\in \M$,  calling  $S(\psi)$ the set of initial points $\bar{x}\in \R^n$ for which the optimization
problem (\ref{oc1})-(\ref{oc2}) has two optimal trajectories ending at distinct terminal points,  we have
\bel{z-S}
\meas\big(S(\psi)\big)~\leq~\limsup_{\nu\to\infty}\big[\meas\bigl(S_\nu(\psi)\bigr)\big]~=~0.
\eeq 
We now observe that, for  every $\bar{x}\in\R^n$, the Pontryagin necessary conditions \cite{BPi,
Cesari, FR} take the form
\bel{PMP1}\left\{
\bega{rl} \dot x&=~f\bigl(x,  u(x,p)\bigr),\\[3mm]
\dot p&=~- p\cdot  f_x\bigl(x,  u(x,p)\bigr)- L_x(x,u(x,p)),\enda
\right.\eeq
with boundary conditions
\bel{PMP3} \left\{ \bega{rl} x(0)&=~\bar x,\\[2mm]
p(T)&=~\nabla \psi\bigl(x(T)\bigr).\enda\right.\eeq
Here the optimal control is determined as the pointwise minimizer
\bel{PMP2} u(x,p)~=~\argmin_{\omega\in \R^m} \Big\{ L(x,\omega) + p\,f(x,\omega)\Big\}.\eeq
By assumptions, $f$ is affine w.r.t.~$\omega$, while by (\ref{u-convex-L})  the cost function $L$ is uniformly convex. 
As a consequence,  the minimizer in (\ref{PMP2})  is unique.  Therefore, the map $(x,p)\mapsto u(x,p)$ is well defined and continuously differentiable, and  the system of ODEs (\ref{PMP1}) has $\C^1$ right hand sides. 
We conclude that, for any $y\in\R^n$, the system (\ref{PMP1}) with  terminal conditions
\bel{td3} x(T)~=~y,\qquad p(T)~=~\nabla\psi(y),\eeq
admits a unique  solution  $t\mapsto (x,p) (t;y)$ on $[0,T]$. In particular, this implies that if  two optimal trajectories starting from $\bar{x}$ have the same terminal point, then then they must coincide for all $t\in [0,T]$. Hence (\ref{z-S}) yields (ii).
\medskip

{\bf 2.} Given $\nu\geq 1$, we now claim that $\M_\nu$ is open in $\C^2_+(\R^n)$. 
Indeed, thanks to the uniform bounds on optimal controls and optimal trajectories proved in Lemma~\ref{prio},
standard arguments show that each set $S_\nu(\psi)$ is closed and bounded.   
Moreover, since the minimum cost
for (\ref{oc1})-(\ref{oc2}) depends continuously on $\bar x, f,L,\psi$, the map $ \psi\mapsto S_\nu(\psi)$
is upper semicontinuous.

Given any terminal cost $\Tilde{\psi}\in \M_\nu$, let $A$ be an open set  such that 
$$S_\nu(\Tilde{\psi})\,\subset\, A,\qquad\qquad \meas(A)<{1\over \nu}\,.$$
Based on Lemma~\ref{bb},  for any initial datum $\bar x\in \ov B_\nu$, every optimal control
$u^*$ and optimal trajectory $x^*$ satisfy
\bel{uxbound}\hbox{\rm ess-}\!\!\!\sup_{t\in [0,T]}  |u^*(t)|~\leq~\alpha(\nu),\qquad\qquad 
\sup_{t\in [0,T]}|x^*(t)|~\leq~\beta(\nu).
\eeq
By upper semicontinuity, there exists $\delta_0>0$ such that 
$$
\|\psi-\Tilde\psi\|_{\C^2} ~<~\delta_0\qquad\Longrightarrow\qquad S_\nu(\psi)\,\subset\, A,\qquad \inf_{x\in\R^n}\Tilde{\psi}>0.$$
As a consequence, $\psi\in \M_\nu$, proving our claim.
\medskip

{\bf 3.}    In the remaining steps, we prove that each $\M_\nu$ is dense in $\C^2_+(\R^n)$. 
Given any $\psi
\in \C^2_+(\R^n)$, we shall construct a small perturbation of $\psi$ that lies inside $\M_\nu$. 

Using Lemma~\ref{bb}, we choose a radius $\rho>0$ large enough
so that the ball $B_{\rho}$  contains all trajectories that satisfy the PMP (\ref{PMP1}) and start  at some point $\bar x\in \ov B_\nu$.  

Denoting by $t\mapsto (x,p) (t;y)$  the  unique  solution of the system of ODEs (\ref{PMP1})
with terminal data (\ref{td3}),  we observe that the map $y\mapsto x(0;y)$ is $\C^1$. 
%
Consider the sets
\bel{Sd0} S_{\delta_0}~\doteq~\Big\{ x(0,y)\,;~~y\in B_{\rho}, \quad \bigl| \det (D_y x(0,y))\bigr| \leq \delta_0\Big\},\quad\qquad S^{-1}_{\delta_0}\,\doteq\,\big\{y\,;~ x(0,y)\in S_{\delta_0}\big\}.
\eeq
By choosing $\delta_0>0$ sufficiently small we obtain
\bel{S-delta}
\meas \bigl( S_{\delta_0}\bigr)~\leq~\delta_0\cdot \meas(B_{\rho})~<~{1\over 2\nu}.
\eeq
Next, consider the open subset of couples in $\R^{n+n}$
\bel{G-nu}
\Gamma_{\nu}~\doteq~\left\{ (y_1, y_2)\in \left(B_{\rho}\backslash S^{-1}_{\delta_0}\right)\times  \left(B_{\rho}\backslash S^{-1}_{\delta_0}\right)\,;~~|y_1-y_2|> {9\over 10\nu}\right\}.
\eeq
%
%
For every couple of points $(\bar y_1,\bar y_2)\in\Gamma_{\nu}$, let $\vp^{(\bar y_1, \bar y_2)}\in\C^\infty_c(\R^{n})$ be a smooth function with compact support
such that 
\bel{vpp} \vp^{(\bar y_1, \bar y_2)} (y) ~=~\left\{ \bega{rl} 1\quad 
&\hbox{if}\quad |y-\bar y_1|\leq {1\over 5\nu},\\[3mm]
-1\quad 
&\hbox{if}\quad |y-\bar y_2|\leq {1\over 5\nu},\\[3mm]
0\quad &\hbox{if}\quad |y-\bar y_1|\geq {2\over 5\nu}~~\hbox{and}~~|y-\bar y_2|\geq {2\over 5\nu}. \enda\right.
\eeq
Covering  the compact closure $\ov\Gamma_{\nu}$ with finitely many balls, say $B\left(\bigl(y^k_1,y^k_2\bigr), {1\over 5\nu}\right)$ for $(y^k_1,y^k_2)\in \ov \Gamma_\nu$, $k\in\{1,\dots, N_{\nu}\}$, we define  a family of terminal costs, depending on the additional parameters $\theta=\bigl(\theta_{1},\dots, \theta_{N_{\nu}}\bigr)$:
\bel{psi-t}
\psi^{\theta}(y)~=~\psi(y)+ \sum_{k=1}^{N_{\nu}} \theta_{k}\cdot\vp^{(y^k_1,y^k_2)}(y)\qquad\forall y\in\R^n.
\eeq

{\bf 4.} 
For any given $\theta\in \R^{N_{\nu}}$, let $(x^{\theta}(s;y),p^{\theta}(s;y))$ be the solution of (\ref{PMP1}) with  terminal condition $(x(T),p(T))=\big(y, \nabla\psi^{\theta}(y)\big)$.   We denote by  $J^{\theta}(y)$  the cost of this trajectory:
\[
J^{\theta}(y)~=~\int_{0}^{T} L\big(x^\theta(t;y), u\big(x^{\theta}(t;y),p^{\theta}(t;y)\big)\big)dt+\psi^{\theta}(y).
\]
Observe that, for any   $k\in\{1,\dots, N_{\nu}\}$ and any  $y\in B\left(y^{k}_1,{1\over 5\nu}\right)\bigcup B\left(y^{k}_2,{1\over 5\nu}\right)$, the definition (\ref{vpp}) implies
\bel{n-phi}
\nabla\psi^{\theta}(y)~=~\nabla\psi(y)+\sum_{j\in \{1,\dots, N_{\nu}\}\backslash\{k\}}\theta_j\cdot \nabla\vp^{(y^j_1,y^j_2)}(y).
\eeq
In this case,  $\left(x^{\theta}(t;y),p^{\theta}(t;y)\right)$ does not depend on $\theta_k$ and 
\bel{d-t-k}
{\partial x^{\theta}\over \partial \theta_k}(0;y)~=~0,\qquad\quad {\partial J^{\theta}\over \partial \theta_k}(y)~=~{\partial \psi^{\theta}\over \partial \theta_k}(y)~=~\begin{cases}~1&\mathrm{if}\quad y\in  B\left(y^{k}_1,{1\over 5\nu}\right),\\[2mm]
-1&\mathrm{if}\quad y\in  B\left(y^{k}_2,{1\over 5\nu}\right).
\end{cases}
\eeq

\medskip
{\bf 5.}
Define the map $\Phi: \Gamma_{\nu}\times \R^{N_{\nu}}\to \R^{n+1}$ by setting
\bel{Phdef}
\Phi(y_1,y_2,\theta)~=~\left(x^{\theta}(0;y_1)-x^{\theta}(0;y_2), J^{\theta}(y_1)-J^{\theta}(y_2)\right).
\eeq
 for all $(y_1,y_2,\theta)\in \Gamma_{\nu}\times \R^{N_{\nu}}$.
For any $k\in \{1,2,\dots, N_{\nu}\}$, by (\ref{d-t-k}) it now follows
\[
{\partial \Phi \over \partial\theta_k}(y_1,y_2,\theta)~=~(0, 0, 2)\qquad\forall (y_1,y_2)\in B\left(\big(y^k_1,y^k_2\big), {1\over 5\nu}\right).
\]
Moreover, by (\ref{Sd0}) and (\ref{psi-t}),  there exists  $\delta_1>0$ small enough such that 
\bel{Det}
\Big|\mathrm{det}\bigl(D_y x^{\theta}(0;y_i)\bigr)\Big|~>~{\delta_0\over 2},\qquad i\in \{1,2\}, ~~\theta=(\theta_1,\ldots,\theta_{N_\nu})\in B_{\delta_1}.
\eeq
Therefore, $\mathrm{rank}\Big( D\Phi(y_1,y_2,\theta)\Big)=n+1$ and 
$\Phi$ is transversal to  the zero manifold
$$\{(0,0)\}~\subset~\R^n\times\R$$ on $B\left((y^k_1, y^k_2),\, {1\over 5\nu}\right)\times B(0,\delta_1)$. Since these balls provide a covering, we conclude that $\Phi$ is transversal to $\{(0,0)\}$ on the whole domain $\Gamma_{\nu}\times B(0,\delta_1)$. 
\medskip

{\bf 6.} Finally, by the transversality theorem \cite{Bloom, GG},  there exists a set $\Theta\subset\R^{N_\nu}$, dense in the ball $B(0,\delta_1)$, 
such that for every $\theta\in \Theta$ the map $\Phi(\cdot,\cdot,\theta)$ is transversal to the zero manifold
$\big\{ (0,0)\}\subset\R^n\times\R$.    
This means: for every couple $(\bar y_1,\bar y_2)\in \Gamma_{\nu}$
such that 
$$x^{\theta}(0,\bar y_1)= x^{\theta}(0, \bar y_2),\qquad\quad J^{\theta}(\bar y_1) = J^{\theta}(\bar y_2),$$
the Jacobian $D_{(y_1,y_2)}\Phi(\bar y_1,\bar y_2,\theta)$ has rank $n+1$.   
Hence, by the implicit function theorem, the set 
of couples 
$$
\Gamma_{(\bar{y}_1,\bar{y}_2)}(r)~\doteq~\Big\{ (y_1,y_2)\in B\big(\big(\bar{y}_1,\bar{y}_2\big),r\big)\cap \Gamma_{\nu}\,;~~\Phi(y_1,y_2,\theta)=(0,0)\in \R^{n+1}\Big\}
$$
is contained in an $(n-1)$-dimensional  manifold, for some $r>0$ small. 
The $n$-dimensional measure of this set is thus
\[
\meas\left(\bigl\{x^\theta(0;y_1)\,;~~(y_1,y_2)\in \Gamma_{(\bar{y}_1,\bar{y}_2)}(r)\bigr\}\right)~=~0.
\]
In turn, for every $\theta\in \Theta$ this implies 
\[
\meas\Big(\bigl\{x^{\theta}(0;y_1)\in \R^n;~~ \hbox{there exists $y_2\in\R^n$ such that}~~(y_1,y_2)\in \Gamma_{\nu}, ~\Phi(y_1,y_2,\theta)=0\bigr\}\Big)~=~0.
\]
On the other hand,  since there exists a constant $C>0$ such that 
\[
\big|x^{\theta}(0,y)-x(0,y)\big|~\leq~C|\theta|\qquad\forall \theta\in\Theta, ~y\in  S^{-1}_{\delta_0},
\]
we have 
\bel{mba}
\left\{ x^\theta(0,y)\,;~~y\in   S^{-1}_{\delta_0}\right\}~\subseteq~B\bigl(S_{\delta_0},\,C|\theta|\bigr).
\eeq
Since $S_{\delta_0}$ is compact, the measure of the $\ve$-neighborhood around the set $S_{\delta_0}$ satisfies
$$\lim_{\ve\to 0}~ \meas 
\Big(B(S_{\delta_0},\ve)\Big)~=~ \meas (S_{\delta_0}).$$
Therefore, choosing $|\theta|$ small enough,  by (\ref{S-delta}) and (\ref{mba}) we obtain
$$\bega{l}\ds\meas\bigl(S_\nu(\psi^{\theta})\bigr)~\leq~\meas\Big(\bigl\{ x^\theta(0,y)\,;~~y\in   S^{-1}_{\delta_0}\bigr\}\Big)\\[3mm]\qquad\ds ~\leq~\meas\left(B\bigl(S_{\delta_0},\,C|\theta|\bigr)\right)
~<~\meas\bigl(S_{\delta_0}\bigr)+{1\over 2\nu}~<~{1\over\nu}\,.\enda$$
Hence the terminal cost $\psi^{\theta}$ lies in $\M_\nu$.   This shows that $\M_\nu$ 
is everywhere dense, completing the proof.
\endproof

\section{Non-uniqueness for  mean field games}
\setcounter{equation}{0}
\label{s:5}

Consider again the optimal control problem (\ref{oc1})-(\ref{oc2}), with $f,L$ satisfying 
{\bf (A1)-(A2)} and $\psi\in \C^2_+$.
Let $Y^*(t)\doteq \begin{pmatrix} x^*(t)\cr p^*(t)\end{pmatrix}$ be a solution to Pontryagin's optimality conditions
(\ref{PMP1})-(\ref{PMP2}). 
 Linearizing the system of ODEs in (\ref{PMP1}) at $Y^*$, we obtain a system of the form
\bel{linY}\dot Y~=~A(t) Y,\eeq
describing the evolution of a first order perturbation.
We shall assume that $Y(t)=\begin{pmatrix} X(t)\cr P(t)\end{pmatrix}\equiv 0$ is the only solution to the linearized system (\ref{linY}) with boundary 
conditions
\bel{lbc} X(0)~=~0,\qquad \qquad P(T)~=~D^2\psi(x^*(T))\cdot X(T).\eeq
Notice that these assumptions imply that this solution is structurally stable.  By the implicit function theorem,
one can slightly perturb the dynamics and the cost function, and still find a solution to the 
equations (\ref{PMP1})--(\ref{PMP2}) close to $Y^*$.

In this setting, it is easy to construct a MFG where $x(t,\xi)=x^*(t)$  is a structurally stable solution.
Indeed, define the barycenter
\bel{bar} b(t)~\doteq~\int_0^1 x(t,\xi)\, d\xi.\eeq
Consider a game where the state of each player evolves  with the same dynamics
\bel{mf3} \dot x(t)~=~f\bigl(x(t), u(t)\bigr)~=~f_0(x(t))+\sum_{i=1}^mf_i(x(t))\cdot u_i(t),\qquad x(\xi,0)= \bar x, \qquad \xi\in [0,1],\eeq
and all players share the same cost functional
\bel{mf4}
J~\doteq~\int_0^T \Big[L\bigl( x(t), u(t)\bigr) +\kappa \bigl|x(t) - b(t)\bigr|^2\Big]\, dt 
+ \psi\bigl(x(T)\bigr).
\eeq
\begin{theorem}\label{t:31} Assume that $f,L$ satisfy
{\bf (A1)-(A2)} while $\psi\in \C^2_+$. Let $(x^*, p^*)$ be a solution to the Pontryagin equations
(\ref{PMP1})-(\ref{PMP2}).
Then, if the constant $\kappa>0$ is large enough,
the MFG (\ref{mf4})-(\ref{mf3}) admits a solution where $x(t,\xi)= x^*(t)$ for all $\xi\in [0,1]$, $t\in [0,T]$.  

If the linearized system (\ref{linY})-(\ref{lbc}) has only the zero solution,  then this solution of the MFG is structurally stable.
\end{theorem}
{\bf Proof.}
{\bf 1.} W.l.o.g., we can assume $\psi=0$. Indeed, the above optimal control problem can always be written as a Bolza problem, replacing the functional $J$ at (\ref{oc1}) with
\bel{OBP} J^\sharp~\doteq~\int_0^T \bigl[L(x,u) + \nabla \psi(x) \cdot f(x,u)\bigr]\, dt.\eeq
If $\bigl(x(t), p(t), u(t)\bigr)$ provide a solution to the equations (\ref{PMP1})-(\ref{PMP2})
for the original problem, one readily checks that the triple
$\bigl(x(t), p(t) - \nabla \psi(x(t)), u(t)\bigr)$ provides a solution to the corresponding Pontryagin's equations
for the Bolza problem (\ref{OBP}).
\v
{\bf 2.} We thus assume that  $\psi=0$. For every given  $b\in \C^{0}([0,T])$ with $\|b-b^*\|_{\C^0}\leq 1$, we claim that (\ref{mf3})-(\ref{mf4}) admits a unique optimal solution for $\kappa>0$ sufficiently large. Indeed, let $(u^b,x^b)$ be a pair of optimal control and  optimal trajectory of (\ref{mf3})-(\ref{mf4}). By  Lemma \ref{prio}, it follows 
\bel{bd-1}
\|u^b\|_{{\bf L}^{\infty}}, \|x^{b}\|_{\C^{0}} ~\leq~C_1,\qquad \kappa \cdot \int_{0}^{T} |x^b(t)-b(t)|^2dt~\leq~C_1,
\eeq
for some  $C_1>0$   which depends only on $f,L$ and $T$. By the necessary conditions, there exists $p^b\in C^{0}([0,T])$ such that $(x^b,p^{b}, u^b)$ solves the PMP
\bel{PMP-non}\left\{
\bega{rl} \dot x&=~f\bigl(x,  u(x,p)\bigr),\\[3mm]
\dot p&=~- p\cdot  f_x\bigl(x,  u(x,p)\bigr)-L_x(x,u(x,p))- 2\kappa(x-b)\,,\enda
\right.
\eeq
with $x(0)=\bar{x}$, $p(T)=0$ and, recalling (\ref{PMP2}),  
\[
u^b(t)~=~u\bigl(x^b(t),p^b(t)\bigr),~\qquad t\in [0,T].
\]
By (\ref{bd-1}) and the second equation of (\ref{PMP-non}) we deduce
\bel{bd-p}
\|p^b\|_{\C^{0}}~\leq~C_2\left(1+\kappa\cdot \int_{0}^{T}|x^b(t)-b(t)|dt\right)~\leq~C_2\left(1+\sqrt{\kappa}\right).
\eeq

{\bf 3.} Next,
consider the Hamiltonian
\bel{Ham1}
H^{b}(x,u,p,t)~\doteq~L(x,u) + \kappa \bigl|x-b(t)\bigr|^2 +p\cdot f(x,u),\eeq
and the reduced Hamiltonian
\bel{Ham2}
\Hat H^{b}(x,p,t)~\doteq~\min_{u\in\R^m}\Big\{ L(x,u) + \kappa \bigl|x-b(t)\bigr|^2 +p\cdot f(x,u)\Big\}.\eeq
The  the optimality condition implies
$$\Hat H^b(x^b, p^b,t) ~=~H(x^b, u^b, p^b,t),\qquad\quad \Hat H^b(x, p,t) ~\geq~H^b(x, u, p,t),$$
\bel{H-x}
 \partial_x\Hat H^b(x,p,t)~=~p\cdot f_x\bigl(x,  u(x,p)\bigr)+L_x(x,u(x,p))+ 2\kappa(x-b(t)),
\eeq
and
\[
L_{uu}(x,u^b(x,p))u^b_{p}(x,p)+(f_1(x),f_2(x),\dots, f_m(x))~=~0.
\]
By the uniform convexity of $L(x,u)$ w.r.t $u$ and the bounds on the vector fields $f_i$, it follows
\[
\|\partial_p u(x,p)\|_{\C^{0}}~\leq~C_3\cdot {1\over \delta_L}.
\]
Therefore, from (\ref{H-x}), (\ref{bd-p}) and (\ref{bd-1}), one obtains
\[
\partial_x\Hat H^b(x^b(t), p^b(t),t)~=~-\dot{p}^b(t),\qquad \partial_{xx}\Hat H^b(x, p^b(t),t)~=~2\kappa I+ G(x,t)
\]
with 
\[
\|G(x,t)\|_{C^0}~\leq~C_4\cdot (1+\sqrt{\kappa})\qquad\forall t\in [0,T], |x|\leq C_1.
\]
In particular, for $\kappa>0$ sufficiently large,  the map $x\mapsto \Hat H^b(x, p^b(t),t)$ is strictly convex in $B(0,C_1)$ for all $t\in [0,T]$ and 
\[
 \partial_{xx}\Hat H^b(x, p^b(t),t)~\geq~\kappa I\,.
\]
In this setting, we show that $ x^b$ is the unique optimal solution of (\ref{mf3})-(\ref{mf4}). 
Indeed, let $(u_1,x_1)$ be another  pair of optimal control and  optimal trajectory for 
(\ref{mf3})-(\ref{mf4}). Notice that $\|x_1-b\|_{\C^0}\leq C_1 $. Using the convexity of $\Hat H$ in the variable $x$,  the difference in costs is estimated by 
$$\bega{l}\ds \int_0^T \Big[ L(x_1,u_1)+\kappa |x_1- b|^2-  L(x^b, u^b) -\kappa |x^b- b|^2 \Big]\, dt \\[4mm]
\ds\quad =~\int_0^T \bigl[  H^b(x_1,u_1,p^b,t) -H(x^b, u^b, p^b,t) \bigr]\, dt
- \int_0^T p^b(t) \cdot \bigl[ f(x_1,u_1) -f(x^b, u^b) \bigr]\, dt\\[4mm]
\ds\quad \geq~\int_0^T \bigl[ \Hat H^b(x_1, p^b,t) - \Hat H^b(x^b,p^b,t)\bigr]\, dt
- \int_0^T p^b(t) \cdot \bigl[ \dot x_1(t)- \dot x^b(t)\bigr]\, dt\\[4mm]
\ds\quad \geq~\int_0^T \partial_x \Hat H^b(x^b,p^b,t)\cdot \bigl(x_1(t)-x^b(t)\bigr), dt
- \int_0^T p^b(t) \cdot \bigl[ \dot x_1(t)- \dot x^b(t)\bigr]\, dt\\[4mm]
\ds\quad =~\int_0^T \Big[- \dot p^b(t)\cdot \bigl(x_1(t)-x^b(t)\bigr) - p^b(t) \cdot \bigl( \dot x_1(t)- \dot x^b(t)\bigr)\
\Big], dt\\[4mm]
\quad =~p^b(0)\bigl( x_1(0)-x^b(0)\bigr) - p^b(T)\bigl(x_1(T)-x^b(T)\bigr)~=~0.
\enda
$$
Notice that if $H$ is strictly convex, then one of the above inequality is strict whenever $x_1(t)\not= x^b(t)$.
In this case, the optimal control is unique.
\medskip

{\bf 4.} By the same argument used in Step {\bf 3}, one can show  that, for $\kappa>0$ sufficiently large, $x^*(\cdot)$ is the unique optimal solution of (\ref{mf3})-(\ref{mf4}) with $b=b^*$. In particular, $x^*=b^*$ and  the corresponding control $u^*(t)=u(x^*(t),p^*(t))$ provide the one and only optimal solution for every player. It remains to show that this solution of the MFG is structurally stable. 

Consider the best reply map $b(\cdot) \mapsto \Phi(b)$, defined by
\[
\Phi(b)(t)~=~\int_{0}^{1}x^b(t,\xi)d\xi\qquad\forall t\in [0,T],
\]
where $x^{b}(\xi,\cdot)$ denotes the unique optimal solution of (\ref{mf3})-(\ref{mf4}). 
In this step we show that the linearization of this map
at $b=b^*$ has eigenvalues all $\not= 1$.  Fix ${\bf b}\in C^{0}[0,T]$ with $\|{\bf b}\|_{\C^0}=1$.  For any $\ve\in\R$  sufficiently small, let  $x^{\ve}(t) $ be the unique optimal solution  (\ref{mf3})-(\ref{mf4}) with $b=b^*+\ve{\bf b}$.  By the necessary conditions, there exists $p^{\ve}\in C^{0}([0,T])$ such that $(x^{\ve},p^{\ve})$ solves  PMP (\ref{PMP-non}) with $b=b^*+\ve{\bf b}$. By a linearization, one obtains
\[
\begin{bmatrix}
x^{\ve}(t)\\[2mm]
p^{\ve}(t)
\end{bmatrix}
~=~\begin{bmatrix}
x^*(t)\\[2mm]
p^*(t)
\end{bmatrix}+\ve \, \begin{bmatrix}
X_{{\bf b}}(t)\\[2mm]
P_{{\bf b}}(t)
\end{bmatrix}+o(\ve).
\]
Here, $Y_{\bf b}(t)=\begin{bmatrix}
X_{{\bf b}}(t)\\[3mm]
P_{{\bf b}}(t)
\end{bmatrix}$ is the solution to  the equation obtained by linearizing (\ref{PMP-non})  around $Y^*$, namely
\[
\dot{Y}(t)~=~A(t)Y(t)+2\kappa\cdot \begin{bmatrix}
0\\[3mm]
X-{\bf b}
\end{bmatrix}
\]
with boundary  conditions (\ref{lbc}). 

Let now  $(\lambda,{\bf b})$ be a pair of eigenvalue and eigenfunction of $D\Phi(b^*)$. We then have 
$$
D\Phi(b)({\bf b})~=~X_{{\bf b}}~=~\lambda{\bf b},
$$
and this implies that $Y_{b}(t)$ solves the linear ODE 
\[
\dot{Y}(t)~=~A(t)Y(t)+2\kappa \left(1-{1\over \lambda}\right)\cdot  \begin{bmatrix}
0\\[3mm]
{\bf b}
\end{bmatrix}.
\]
Thus, by the assumption at (\ref{linY})-(\ref{lbc}), it follows $\lambda\not=1$.
\v
{\bf 5.}  To prove the structural stability of the solution to the MFG, 
for $\delta>0$ sufficiently small we consider the perturbed problem 
\bel{ga1-p-0}
\hbox{minimize:} \quad \int_0^T \left[L\bigl( x(t), u(t)\bigr)+  \kappa\cdot
\bigl| x(t) - b(t)\bigr|^2+\delta\cdot L_1(x,u,b)\right]\, dt,\eeq
subject to 
\bel{ga2-p-0}\dot x(t) ~=~f(x,u)+\delta\, g(x,u,b),\qquad\quad x(0)\,=\bar{x}+\, \delta \, \bar x_1(\xi).\eeq
We here assume
\bel{fLx-0}
\|g\|_{\mathcal{C}^2}+\|L_1\|_{\mathcal{C}^2}~\leq~ 1,\qquad \|\bar x_1\|_{{\bf L}^{\infty}}~\leq~1.
\eeq

By the same argument used in Step {\bf 3}, for $\kappa>0$ sufficiently large, the optimal control problem (\ref{ga1-p-0})-(\ref{ga2-p-0})
admits a unique solution for all $\delta\in [0,1]$ and $\|b-b^*\|_{\C^0}\leq 1$.  Consider the best reply map $b(\cdot) \mapsto \Phi^{\delta}(b)$, defined by 
\[
\Phi^{\delta}(b)(t)~=~\int_{0}^{1}x_{\delta}^{b}(t,\xi)d\xi\qquad\forall t\in [0,T],
\]
where $x_{\delta}^b(\xi,\cdot)$ denotes the unique optimal solution of (\ref{ga1-p-0})-(\ref{ga2-p-0}). We 
claim that  there exists a constant  $C_6>0$, independent of $b$ and $\delta$, such that
\bel{Phid-0}
\big\|\Phi^\delta(b)-\Phi(b)\big\|_{\C^0}~\leq~\left\|\int_0^1x^{b}_{\delta}(\cdot,\xi)\,d\xi-x^b(\cdot)\right\|_{\C^0}~\leq~C_6\cdot \delta^{2/3}.
\eeq

 Calling $y^b(\xi,\cdot )$ and $y^b_{\delta}(\xi,\cdot )$ the solution to (\ref{ga2-p-0}) corresponding to $u=u^b$ and  the 
optimal control $u=u^{b}_{\delta}$ respectively but with initial data $x(0)=\bar{x}$,  we have 
\bel{df-1}
\big\|y^b(\xi,\cdot )-x^b(\xi,\cdot )\big\|_{\C_0} ,\big\|y^b_{\delta}(\xi,\cdot )-x^b_{\delta}(\xi,\cdot )\big\|_{\C_0}~\leq~\O(1)\cdot \delta.
\eeq
Since $\left(x_{\delta}^b(\xi,\cdot), u^{b}_{\delta}(\xi,\cdot)\right)$ is the optimal pair of (\ref{ga1-p-0})-(\ref{ga2-p-0}), one has
\begin{equation*}
\begin{split}
2\delta T&~\geq~\delta\cdot \int_{0}^{T} L_1\big(y^b(\xi,t),u^b(t),b(t)\big)-L_1\bigl( x^b_{\delta}(t,\xi), u^{b}_{\delta}(t,\xi),b(t)\bigr)dt\\
&~\geq \ds \int_0^T L\bigl( x^b_{\delta}(t,\xi), u^{b}_{\delta}(t,\xi)\bigr)+  \kappa\cdot
\bigl| x^{b}_{\delta}(t,\xi) - b(t)\bigr|^2-L\bigl( y^b(t), u^{b}(t)\bigr)-  \kappa\cdot
\bigl| y^{b}(t) - b(t)\bigr|^2\, dt,
\end{split}
\end{equation*}
and (\ref{df-1}) implies  
$$
\O(1)(1+\kappa)T\delta~\geq~\int_0^T L\bigl( y^\delta(t,\xi), u^{b}_{\delta}(t,\xi)\bigr)+  \kappa\cdot
\bigl| y^\delta(t,\xi) - b(t)\bigr|^2-L\bigl( x^b(t), u^{b}(t)\bigr)-  \kappa\cdot
\bigl| x^{b}(t) - b(t)\bigr|^2\, dt.
$$
Following the same argument in Step 3, we estimate 
\begin{equation*}
\begin{split}
\O(1)(1+\kappa)T\delta&~\geq~\int_0^T \bigl[ \Hat H(y^\delta(t,\xi), p^b,t) - \Hat H(x^b,p^b,t)\bigr]\, dt
- \int_0^T p^b(t) \cdot \bigl[ \dot y^\delta(t,\xi)- \dot x^b(t)\bigr]\, \\
&~\geq~ \kappa\cdot\int_{0}^{T}\left|y^\delta(t,\xi)-x_b(t)\right |^2dt\,.
\end{split}
\end{equation*}
This yields 
\[
\int_{0}^{T}\left|y^\delta(t,\xi)-x_b(t)\right |^2dt~\leq~\O(1)\cdot \left(T+{T\over \kappa}\right)\cdot \delta.
\]
Notice that 
\[
\|\dot y^\delta(\xi,\cdot)\|_{\C^0},~~ \|\dot{x}_b(\cdot)\|_{\C^0}~\leq~C_5
\]
for some  constant $C_5>0$   which depends only on $f,L$ and $T$. We then have 
\[
\|y^\delta(\xi,\cdot)-x_b(\cdot)\|_{\C^0}~\leq~\O(1)\cdot \delta^{2/3},
\]
and (\ref{df-1}) yields (\ref{Phid-0}).

\medskip

{\bf 6.} We are now ready to complete the proof.  By step {\bf 4}, the eigenvalues $\lambda_n$ of  the compact operator $D\Phi(b^*)$ satisfy
\[
\inf_{n\geq 1} |\lambda_n-1|~\geq~\delta_0~>~0.
\]
As a consequence, the inverse  linear operator $\left[D\Phi(b^*)-{\bf I}\right]^{-1}$ is bounded.  We can define the continuous operator $F$ on $\C^0\bigl([0,T]\bigr)$ as the composition
\[
F(b)~\doteq~b^*+\left[D\Phi(b^*)-{\bf I}\right]^{-1}\circ \left[\Phi(b^*)+D\Phi(b^*)(b-b^*)-\Phi^\delta(b)\right].
\]
From (\ref{Phid-0}) it follows
\[
\|F(b)-b^*\|_{\C^0}~\leq ~{\O(1)\over \delta_0}\cdot \left(\delta+\delta^{2/3}+\|b-b^*\|^2_{\C^0}\right).
\]
Therefore,  for $\delta>0$ sufficiently small, one has
\[
\bigl\|F(b)-b^*\bigr\|_{\C^0}~\leq~\delta^{1/3}\qquad\forall b\in \mathcal{C}\bigl([0,T]\bigr),~~ 
\|b-b^*\|_{\C_0}\leq \delta^{1/3}.
\]
On the other hand, for every $b\in\mathcal{C}([0,T])$ with $ \|b-b^*\|_{\C^0}\leq \delta^{1/3}$, the function $F(b)(\cdot)$ is Lipschitz continuous 
with some uniform Lipschitz constant $M$. In particular, $F$ maps the convex and compact subset \[
 K~=~\left\{b\in\mathcal{C}^0[0,T]: \|b-b^*\|_{\C^0}\leq \delta^{1/3}, ~~\hbox{Lip}(b)\leq M\right\}
 \]
 into itself.
By  Schauder's fixed point theorem, there exists $b_{\delta}\in K $ such that $F(b_{\delta})=b_{\delta}$.  This implies that $b_{\delta}$ is a fixed point of $\Phi^{\delta}$ with $\|b_{\delta}-b^*\|_{\C_0}\leq \sqrt{\delta}$. The family of optimal $x_{\delta}^{b_{\delta}}(\cdot,\xi)$, $\xi\in\Omega$  thus provide a solution to the perturbed MFG, such that 
\[
\sup_{t\in [0,T]}\int_{0}^{t}\big|x_{\delta}^{b^{\delta}}(\cdot,\xi)\big|~d\xi~\leq~T\cdot \delta^{1/3}.
\]
Therefore, the  solution $x(t,\xi)\equiv x^*(t)$ is is structurally stable.
\endproof


%
%
%

\begin{remark}\label{r:31}{\rm
An immediate consequence of the above results is the non-uniqueness of solutions to mean field games.
Namely, given $(f, L)$ satisfying {\bf (A1)-(A2)}, let $\psi\in   \C^2_+(\R^n)$ determine an optimal control problem where, for some $\bar x\in \R^n$, the 
system (\ref{PMP1})-(\ref{PMP2}) admits two distinct solutions, both satisfying the 
structural stability assumptions in Theorem~\ref{t:31}. Then, by choosing $\kappa>0$ large enough,
we obtain a MFG with two solutions, both structurally stable.
In particular, non-uniqueness holds on an open set of MFG.}
\end{remark}

\section{Generic single-valuedness of the best-reply map}
\setcounter{equation}{0}
\label{s:7}
In general, for a given $\eta(\cdot)$ in (\ref{etai}), there will be several players $\xi\in \Omega$
for which
 the optimal control problem (\ref{mmin})--(\ref{iitt})
has multiple solutions. For this reason, the map $\eta\mapsto \tilde\eta = \Phi(\eta)$ at (\ref{mot})
can be multivalued.
Lacking convexity, one cannot guarantee the existence of a fixed point.  
The main result proved in this section is that, for a generic MFG, 
for every $\eta(\cdot)$ in a suitable bounded subset of $\C^2$ functions, 
the set of players $\xi\in \Omega$ having multiple optimal controls has measure zero.
Hence the best reply map (\ref{mot}) is single-valued. The existence of a  fixed point, and the existence of a strong solution to the MFG,  thus 
follow directly from  Schauder's theorem.  
Throughout this section, we consider a quadruple $(\bar x, f,L,\phi)\in \L^\infty\times\C^3\times \C^{3}\times\C^3$ such that  $f,L$, and $\bar x$ satisfy the following assumptions
\begi
\item [{\bf (B1)}] {\it  The function $f$ is affine w.r.t.~the control:
\bel{faff} f(x,u,\eta)~=~f_0(x,\eta)+\sum_{i=1}^m f_i(x,\eta) u_i\,.\eeq
For some constant $c_1$ independent of $\eta$, the vector fields $f_i$ satisfy
\bel{fibo} \bigl|f_i(x,\eta)\bigr|~\leq~c_1\, \bigl(|x|+1\bigr).\eeq
}
\item[{\bf (B2)}] {\it 
There exist constant $c_2>0$ and a continuous function $\ell$ independent of $\eta$  such that, for all $(x,u,\eta)\in \R^n\times \R^m\times \R^N$, one has
\[
\begin{cases}
L(x,u,\eta)~\geq~c_2\bigl( |u|^2-1\bigr),\\[4mm]
\bigl|L_x(x,u,\eta)\bigr|~\leq~\ell(|x|)\cdot (1+|u|^2).
\end{cases}
\]
Moreover, for every $x,\eta$, the map $u\mapsto L(x,u,\eta)$ is uniformly convex.
Namely, for some $\delta_L>0$, 
the $m\times m$ matrix of second derivatives w.r.t.~$u$ satisfies
\bel{Luu}
L_{uu}(x,u,\eta)~>~\delta\cdot \mathbb{I}_m \,,\eeq
for some $\delta>0$, uniformly positive for $x,u,\eta$ in bounded sets.
}
\item[{\bf (B3)}] {\it 
The initial distribution of players, i.e.~the push-forward of the Lebesgue measure on $[0,1]$
via the map $\xi\mapsto \bar x(\xi)\in \R^n$, is a probability measure $\mu_0$ with bounded support and uniformly bounded density w.r.t.~Lebesgue measure on $\R^n$.
}
\endi

Under the above assumptions, by Lemma~\ref{bb} every optimal control $u^*(\cdot)$ and  optimal trajectory $x^*
(\cdot)$ for  the optimization
problem (\ref{mmin})--(\ref{iitt}) satisfy the bounds 
\bel{pri1}
\hbox{\rm ess-}\!\!\!\sup_{t\in [0,T]}  |u^*(t)|~\leq~\alpha_0\doteq \alpha\bigl(\|\bar x\|_{\L^{\infty}}\bigr),\qquad\qquad 
\sup_{t\in [0,T]}|x^*(t)|~\leq~\beta_0\doteq\beta\bigl(\|\bar x\|_{\L^{\infty}}\bigr).
\eeq
As a consequence, any statistic
$\eta(\cdot)$ in (\ref{etai}) will satisfy the a priori bound
\bel{eta-0}
\|\eta\|_{\C^0}~\leq~\gamma_0~\doteq~\left(\sum_{i=1}^{N}\left|\max_{t\in [0,T], |x|\leq \beta(\|\bar x\|_{\L^{\infty}}) }\phi_i(t,x)\right|^2\right)^{1/2}.
\eeq
Next, we recall that,
for any given $\eta(\cdot)$,
by the optimality conditions there exists an adjoint vector
$p^*\in C^{0}([0,T])$ such that $(x^*,p^{*})=(x^{\eta},p^{\eta})$ solves the PMP
\bel{PMP6} \left\{ \bega{rl} \dot x&=~f(x,u^{\eta}(t,x,p),\eta),\\[2mm]
\dot p&=~- p \cdot f_x(x,u^{\eta}(t,x,p),\eta) - L_x\bigl(x(t), u^{\eta}(t,x,p),\eta\big),\enda\right.
\eeq
with terminal data of the form
\bel{PMP4}
x(T)\,=\,y,\qquad\qquad p(T)\,=\, \nabla\psi(y).\eeq
Here the optimal control $u^*(t)=u^{\eta}(t,x^*,p^*)$ is given by
$$u^\eta(t,x,p)~=~\argmin_\omega\Big\{ L\bigl(x,\omega,\eta(t)\bigr) +  p \,f\bigl(x,\omega,\eta(t)
\bigr)\Big\}.$$
By the strict convexity of $L$, since $f$ is affine w.r.t.~$u$, this minimizer can be determined as the unique solution to
\bel{u-eta}
 L_u\bigl(x,\omega,\eta(t)\bigr)+ f_u\bigl(x,\omega,\eta(t)\bigr)~=~0.
\eeq
Relying on the uniform bound on all optimal controls and optimal trajectories, proved in 
Lemma~\ref{bb}, we now establish a uniform bound on all statistics $\eta(\cdot)$ in 
(\ref{etai}).
\begin{lemma}\label{pb-eta}
Under the assumptions {\bf (B1)-(B2)}, for any $\phi\in \C^3$ and any terminal cost $\psi\in \C^2_+$, there exists a constant $\gamma_3$ such that the  composed map $\Phi$ in (\ref{Phi})-(\ref{mot}) satisfies the implication
\bel{phi3}\|\eta\|_{\C^3}\leq \gamma_3\qquad\implies\qquad  \bigl\|\Phi(\eta)\bigr\|_{\C^3}\leq \gamma_3\,.\eeq
\end{lemma}

{\bf Proof.}  {\bf 1.}
By the optimality conditions (\ref{pri1})-(\ref{PMP4}), the adjoint vector $p^*$ is bounded by 
\bel{peta-3}
\|p^*\|_{\C^0}~\leq~\left(\|\nabla\psi\|_{\C^0(B_{\beta_0)}}+\|D_xL\|_{\C^0(B_{\alpha_0+\beta_0+\gamma_0})}\right)\cdot \exp\left(T\cdot\|D_xf\|_{\C^0(B_{\alpha_0+\beta_0+\gamma_0})}\right)\doteq \sigma_0.
\eeq
Therefore, setting $r_0\doteq \alpha_0+\beta_0+\gamma_0+\sigma_0$, 
we can assume that $(x,u,p,\eta)$ take values  inside a fixed ball $B_{r_0}$.
\medskip

{\bf 2.} 
Proceeding by induction, we will show the implications
\bel{ebo}  \|\eta\|_{\C^k}\,\leq\,\gamma_{k}\qquad\implies\qquad  
 \bigl\|\Phi(\eta)\bigr\|_{\C^{k+1}}\,\leq\,\gamma_{k+1}\qquad\hbox{for}~  k=0,\dots, 2\,,\eeq 
for some suitable constants $\gamma_k$.

Indeed, assume that  $\|\eta\|_{\C^k}\leq\gamma_{k}$ for some $k\leq 2$. 
Since   $f,L_u\in\C^3$ and 
$$L_{uu}(x,u,\eta)~\geq~\delta_{r_0}\cdot \mathbb{I}_m$$
for all $(x,u,\eta)\in B_{r_0}$,
the implicit function theorem implies that the solution  $u(x,p,\eta)$ of (\ref{u-eta}) is in $\C^k$ and 
satisfies
\[
\|u\|_{\C^{k}(B_{r_0})}~\leq~\left({1\over \delta_{r_0}}\right)^{km}\cdot \alpha_k\,.
\]
Here the constant  $\alpha_{k}>0$ depends on $r_0$, $\|f\|_{\C^{k+1}(B_{r_0})}$,  $\|L\|_{\C^{k+1}(B_{r_0})}$. Hence, the solution $(x^{\eta},p^{\eta})$ of  (\ref{PMP6}) is in $\C^{k+1}$ and 
\[
\|x^{\eta}\|_{\C^{k+1}}~\leq~\beta_{k+1},\qquad \|p^{\eta}\|_{\C^{k+1}}~\leq~\sigma_{k+1}
\]
with $\beta_{k+1}, \sigma_{k+1}>0$ depending on $r_0$, $\|f\|_{\C^{k+1}(B_{r_0})}$,  $\|L\|_{\C^{k+1}(B_{r_0})}$, and $\gamma_0,\dots,\gamma_k$. As a consequence, 
(\ref{Phi}) implies
\begin{eqnarray*}
\|\Phi(\eta)\|_{\C^{k+1}}~\leq~\gamma_{k+1}
\end{eqnarray*}
where the constant $\gamma_{k+1}>0$  can be computed in terms of   $r_0$, $\|f\|_{\C^{k+1}(B_{r_0})}$,  $\|L\|_{\C^{k+1}(B_{r_0})}$,   $\|\phi\|_{\C^{k+1}(B_{r_0})}$, and $\gamma_0,\dots,\gamma_k$. 

Thus,  by induction, (\ref{eta-0})  yields an a priori bound of $\eta$ in (\ref{ebo}).
In particular, (\ref{phi3}) holds.
\endproof

We are now ready to prove the main result of the paper.

\begin{theorem}\label{t:41} Consider the mean field game at (\ref{etai})--(\ref{iitt}). 
Assume that $(\bar{x},f,L,\phi)\in \L^{\infty}\times \C^3\times \C^{3}\times \C^3$, with 
$f,L,\bar{x}$ satisfying {\bf (B1)-(B3)}.
Then, for any constant $K>0$,  there exists  a  $\Gd$ set $\M\subset  \C^2_+(\R^n)$ such that for every terminal cost $\psi\in\M$, the map $\eta\mapsto \Phi(\eta)$ at (\ref{Phi}) is single-valued on the ball
\bel{Baga}B_K~\doteq~\bigl\{ \eta:[0,T]\mapsto\R^N\,;~~\|\eta\|_{\C^3}\leq K\bigr\}.\eeq
As a consequence, the MFG admits a strong solution.
\end{theorem}

{\bf Proof.}  By suitably choosing the family $\M$ of terminal costs, we need to show that,  if $\psi\in\M$ and $\|\eta\|_{ \C^3}\leq K$, then the set of 
players 
$$P^{\eta}~=~\Big\{\xi\in \Omega\,;~\hbox{the optimal control problem (\ref{mmin})--(\ref{iitt})
has multiple solutions}\Big\}$$
has zero measure. 

Toward this goal, let  $\mu_0$ be a probability measure on $\R^n$ with bounded support and whose density w.r.t.~Lebesgue measure
is uniformly bounded. Assume that for every given  $\ve_0>0$, we can prove
\begi
\item[{\bf (G)}] {\it  There exists  an open dense subset   $\mathcal{M}_{\ve_0}\subset \C^2_+(\R^n)$ such that for every $\psi\in\M_{\ve_0}$ and $\eta\in B_{K}$, 
the set of initial points
\bel{seta} \bega{rl}S^\eta_{\ve_0}&\doteq~\Big\{x_0\in \R^n\,;~~\hbox{the optimization problem (\ref{mmin})--(\ref{iitt}) has two solutions
$x_1(\cdot), x_2(\cdot)$}\\[2mm]
& \qquad\qquad \qquad \hbox{with}~~
 x_1(0)=x_2(0)= x_0\,,\qquad |x_1(T)-x_2(T)|\geq\ve_0\Big\}\enda
 \eeq
has measure
\bel{small}\mu_0(S^\eta_{\ve_0})~<~\ve_0\,.\eeq
}
\endi
Then the set $\M=\ds\bigcap_{\ve_0>0}\M_{\ve_0}$ is  $\Gd$ subset of $\C^2_+(\R^n)$.
Moreover,  for every $\psi\in\M_{\ve_0}$ and $\eta\in B_{K}$ one has
\[
\meas(P^{\eta})~\leq~\lim_{\ve_0\to 0+}\mu_0\left(S^\eta_{\ve_0}\right)~=~0\,.
\]
Indeed, this follows from the observation that, if two optimal trajectories have the same 
terminal point, then by the necessary conditions they must coincide for all $t\in [0,T]$.

In the next several steps, we thus focus on a proof of {\bf (G)}.

\medskip

{\bf 1.}  Given $\ve_0>0$, we claim that the set 
\bel{M-psi}
\mathcal{M}_{\ve_0}~\doteq~\left\{\psi\in\mathcal{C}^2_+(\R^n)\,;~~\mu_0\big(S^{\eta}_{\ve_0}\big)< \ve_0~~ \mathrm{for~every}~\eta\in B_{K}\right\}
\eeq is open.   Equivalently, its complement $\mathcal{M}^c_{\ve_0}$  is closed. 

 Indeed, consider any sequence of elements $\psi_n\in  \mathcal{M}^c_{\ve_0}$ converging to  $\ov \psi$ in $\C^2$ as $n\to\infty$.
For each $n\geq 1$, let $\eta_n\in B_K\subset\mathcal{C}^3$ be such that 
\bel{eta-n}
\mu_0\big(S^{\eta_n}_{\ve_0}\big)~\geq~ \ve_0.
\eeq
By possibly taking a subsequence, we can assume
 that $\eta_n$ converges to $\ov\eta$ in $\mathcal{C}^2$. By the upper semicontinuity of the set of optimal solutions, one has 
\bel{inc1}
\mathcal{S}^{\ov\eta}_{\ve_0}~\supseteq~\limsup_{n\to\infty}\mathcal{S}^{\eta_n}_{\ve_0}~\doteq~\bigcap_{n\geq 1}\bigcup_{k\geq n}\mathcal{S}^{\eta_k}_{\ve_0}.
\eeq
Therefore
\[
\mu_0\big(\mathcal{S}^{\ov \eta}_{\ve_0}\big)~\geq~\mu_0\big(\limsup_{n\to\infty}\mathcal{S}^{\eta_n}_{\ve_0}\big)~\geq~\limsup_{n\to\infty}\mu\big(\mathcal{S}^{\eta_n}_{\ve_0}\big)~\geq~\ve_0,
\]
and this yields $\ov \psi\in \mathcal{M}^c_{\ve_0}$. 
\medskip

{\bf 2.} We will establish the density of the set $\mathcal{M}_{\ve_0}$ in $\mathcal{C}^2$ by constructing 
smooth perturbations of the terminal cost $\psi$ which are very small in the $\C^2$ norm, 
but possibly large in $\C^3$. More precisely, let $\rho_0>0$ be an upper bound for the density of the probability measure $\mu_0$ w.r.t.~Lebesgue measure on $\R^n$. 
Choose a radius $r_0>\|\bar x\|_{\L^\infty}$, so that 
\bel{rr}\hbox{Supp}(\mu_0)~\subset~B(0,r_0).
\eeq
Then choose $R_0>0$ large enough so that, for every $\eta\in B_K$,
every optimal solution starting at a point $x_0\in B(0, r_0)$ remain inside the cube $[-R_0,R_0]^n$. 

Dividing $[-R_0,R_0]^n$ into $\nu= \left(\lfloor {2R_0\over \ve_0}\rfloor+1\right)^n$ smaller cubes with side smaller than $\ve_0$, say $\Gamma_1,\ldots, \Gamma_\nu$, the perturbed terminal cost $\psi^\sharp$ will be defined separately on each cube $\Gamma_k$, so that  the following proper
ties hold.
\begi
\item[(i)] $\psi^\sharp$ coincides with $\psi$ on a neighborhood of the boundary $\partial \Gamma_k$.
\item[(ii)]  For every $k=1,2,\ldots,\nu$ one has the bound
\bel{c2s}
\|\psi^\sharp-\psi\|_{\C^2(\Gamma_k)}~<~\ve_0\,.\eeq
\item[(iii)] There exists an open subset $\Gamma'_k\subset\Gamma_k$ such that
\bel{gpri}
\meas(\Gamma_k\setminus\Gamma'_k)~<~{\ve_1\over \nu}\,,\eeq
\bel{td1}
M_k~<~\Bigl| D^3 \psi^\sharp(x)\bigr|~<~2M_k\qquad\forall x\in \Gamma'_k\,,\eeq
\bel{td2}\Bigl| D^3 \psi^\sharp(x)\bigr|~<~2M_k\qquad\forall x\in \Gamma_k\,.\eeq
\endi
It is clear that, given $\ve_0, \ve_1, M_k$, a function $\psi^\sharp$ with the above properties does exist. Moreover, the  increasing sequence of numbers $M_{k+1}$ will be inductively defined in Step 5 so that $M_{k+1}$ is much larger than $M_k$. 
\medskip

{\bf 3.}   For any given $\eta(\cdot)\in B_K$, we consider the map
\bel{yxm}
y~\mapsto~x^\eta(0,y),\eeq
where $t\mapsto \bigl(x^\eta(t,y),p^\eta(t,y)\bigr)$ is the solution of (\ref{PMP6})
with terminal data (\ref{PMP4}). 
By the assumption {\bf (B3)} on the absolute continuity of the measure $\mu_0$ (describing the initial distribution of players) w.r.t.~Lebesgue measure, 
we can choose $\delta_1>0$
such that the following holds.
Calling $D_y x^\eta(0,y)$ the Jacobian matrix of the map (\ref{yxm}), one has
\bel{jsm}
\mu_0\Big(\left\{ y\in [-R_0,R_0]^n\,;~~\bigl| \det \big(D_y x^\eta(0,y)\big)\bigr| \leq \delta_1\right\}\Big)
~<~{\ve_0\over 2}\,.\eeq

From now on, we shall thus  focus on the  set of points $y\in \R^n$ where 
$\bigl| \det \big(D_y x^\eta(0,y)\big)\bigr| > \delta_1$, so that the map $y\mapsto x^\eta(0,y)$
is locally invertible.

\medskip

\begin{figure}[ht]
\centerline{\hbox{\includegraphics[width=10cm]{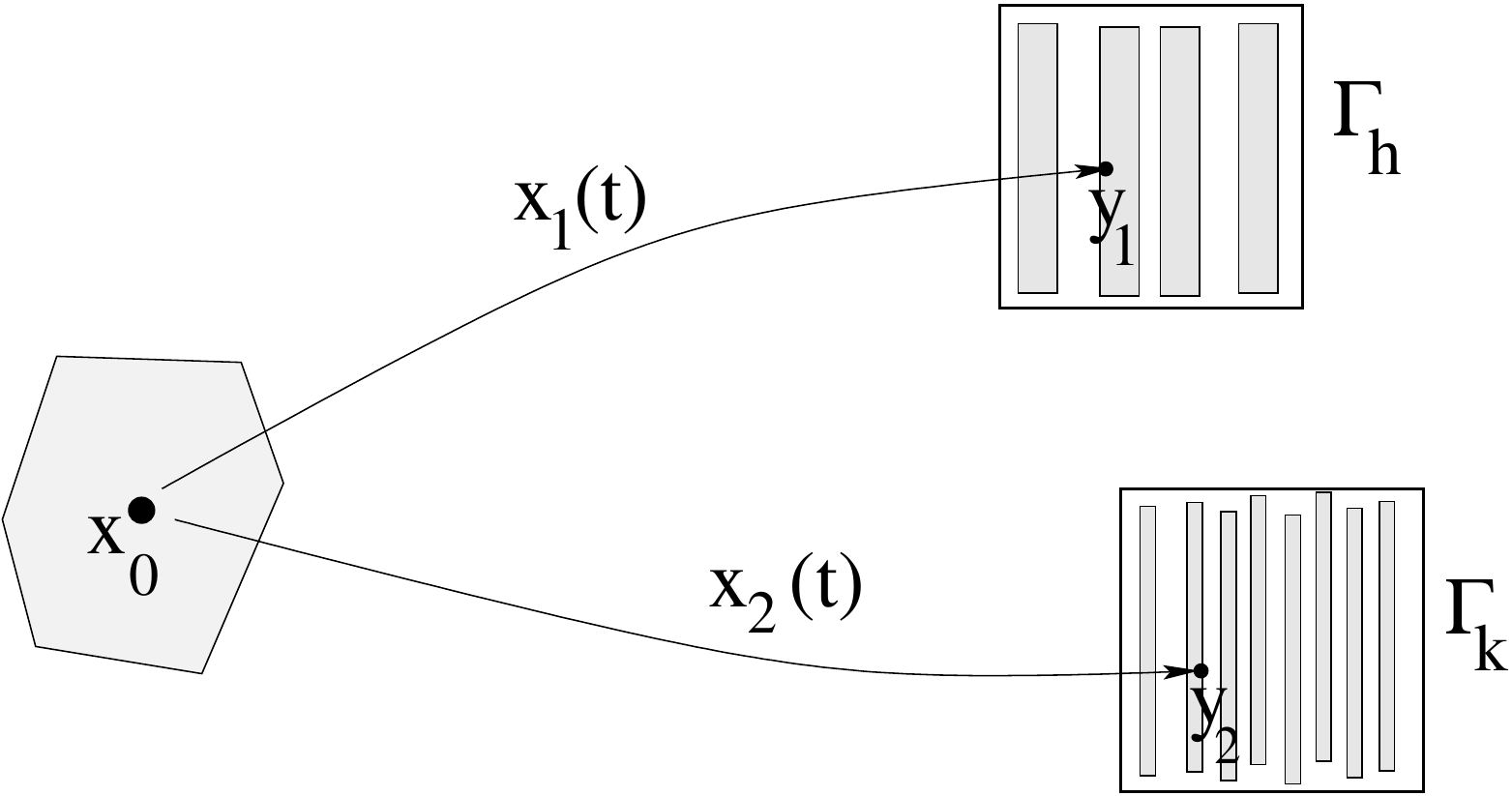}}}
\caption{\small The terminal cost $\psi$ has uniformly bounded gradient.
However, we can construct a perturbation $\psi^\sharp$ whose third derivatives
have vastly different sizes on different cubes $\Gamma_k$ of the partition. 
}
\label{f:mfg2}
\end{figure}
\medskip

{\bf 4.} To help the reader, we first explain the heart of the matter, with the aid of Fig.~\ref{f:mfg2}.     
Let $\eta\in B_K$ be given.
Assume that $x_0$ is an initial point from which two optimal trajectories
$x_1(\cdot)$, $x_2(\cdot)$ 
originate.   To fix ideas, assume
$$y_1\,=\,x_1(T)\in \Gamma'_h\,,\qquad\qquad y_2\,=\,x_2(T)\,\in\, \Gamma'_k\,,$$
with $h<k$. On $\Gamma_k'$ the terminal cost function $\psi^{\sharp}$ has a much larger third
derivative than on
$\Gamma_h'$.  
We observe that the Jacobian matrix of 
the map
$y\mapsto x^\eta(0,y)$ is uniformly invertible in a neighborhood of $y_1$ and $y_2$.
By the implicit function theorem, for all $x\in B(x_0, \delta_2)$, on a ball  centered at $x_0$ with sufficiently 
small radius $\delta_2>0$, we can thus define the cost functions
$\Phi_1(x)$, $\Phi_2(x)$, corresponding to trajectories $x_1(\cdot), x_2(\cdot)$ that start at $x$, satisfy the PMP, and terminate
in a neighborhood of $y_1$, $y_2$, respectively. 
Since the terminal costs $\psi^\sharp\bigl(x_i(T)\bigr)$ of these trajectories have very different third order derivatives, 
we will show that the cost functions $\Phi_1,\Phi_2$ also have different third order 
derivatives in a neighborhood of $x_0$.
Therefore, the set of points where $\Phi_1(x)=\Phi_2(x)$ must be very small, regardless of the particular function $\eta(\cdot)$. 
A proof of these claims will be worked out with the aid of

\begin{lemma}\label{l:42} Consider a system of $n+n$ ODEs on the interval $[0,T]$,
\bel{ODE3}
\left\{\bega{rl} \dot x(t)&=~F\bigl(t,x(t),p(t)\bigr),\\[2mm]
\dot p(t)&=~G\bigl(t, x(t),p(t)\bigr).\enda\right.\eeq
Assume that all coefficients are uniformly bounded in $\C^2$.
\begi
\item[(i)] Consider a family of solutions $(x,p)(t,y)$ with initial data 
\bel{id3}
x(0)= y,\qquad p(0)= \vp(y).\eeq
Assume that  $\vp\in \C^2$ and the map $y\mapsto x(T,y)$ 
is uniformly invertible. More precisely, the norm of its $n\times n$ Jacobian matrix satisfies
\bel{inv} \Big|D_y x(T,y)\Big|\,\leq \, C,\qquad\qquad \Big|\bigl[D_y x(T,y)\bigr]^{-1}\Big|\,\leq \, C.
\eeq
Then the second derivatives  $D^2_x p$ of the map $x(T,y)\mapsto p(T,y)$ satisfy a uniform bound,
depending on the $\C^2$ norms of the functions $F,G, \vp$, and on the constant $C$ in (\ref{inv}).
\item[(ii)] Similarly, consider a family of solutions $(x,p)(t,y)$ with terminal data 
\bel{id4}
x(T)= y,\qquad p(T)= \vp(y).\eeq
Assume that  $\vp\in \C^2$ and the map $y\mapsto x(0,y)$ 
is uniformly invertible. More precisely, the norm of its $n\times n$ Jacobian matrix satisfies
\bel{inv2} \Big|D_y x(0,y)\Big|\,\leq \, C,\qquad\qquad \Big|\bigl[D_y x(0,y)\bigr]^{-1}\Big|\,\leq \, C.
\eeq
Then the second derivatives $D^2_x p$ of the map $x(0,y)\mapsto p(0,y)$ satisfy a uniform bound,
depending on the $\C^2$ norms of the functions $F,G, \vp$, and on the constant $C$ in (\ref{inv2}). 
\endi
\end{lemma}

{\bf Proof.} Part (ii) is entirely similar to part (i), after reversing the direction of time.
We thus focus on a proof of (i).   

Standard results on  the higher order differentiability of 
solutions to ODEs, see for example Theorem 4.1 in \cite{H}, p.100, imply that the maps
\bel{yxp}y\mapsto x(T,y),\qquad\qquad y\mapsto p(T,y)\eeq
are twice continuously differentiable, and satisfy bounds of the form
$$\bigl|D^2_y x(T,y)\bigr|~\leq~C_1,\qquad\qquad \bigl|D^2_y p(T,y)\bigr|~\leq~C_1,
$$ 
for some constant $C_1$ depending only on the $\C^2$ norms of $F,G,\vp$.
By assumption, the first map in (\ref{yxp}) is invertible because of (\ref{inv}).
As a consequence, the inverse function $x\mapsto y(x)$ is well defined,
and has a bounded  second derivatives, depending on the constants $C, C_1$.

This implies that the composed map $x\mapsto p\bigl(T, y(x)\bigr)$
is $\C^2$, and its second derivatives can be bounded in terms of the constants $C, C_1$.
\endproof

We now resume the proof of Theorem~\ref{t:41}.
\v
{\bf 5.} We finalize the construction of the perturbed terminal cost $\psi^\sharp$ by assigning the  increasing sequence of numbers $M_k$.

%
%

We start by choosing $M_1 > \|\psi\|_{\C^3}$.  By induction, assume now that 
$M_1,\ldots, M_{k-1}$ have been chosen.   

Consider any trajectory satisfying the PMP, starting at some point $x\in B_R$ and ending
inside some $\Gamma_j$ with $j\leq k-1$.
We shall apply  Lemma~\ref{l:42} in the special case where (\ref{ODE3}) is given 
by (\ref{PMP6}).

Calling  $t\mapsto (x_j,p_j)(t, y)$ the solution to (\ref{PMP6}) with terminal condition
$(x(T),p(T))=(y,\psi^{\sharp}(y))$ for $y\in \Gamma_j$.  For $x=x_j(0, y)$, we define
\[
\Phi_j(x)~\doteq~\int_{0}^{T}L\Big(x_j(t,y),u^\eta(t,x_j(t,y),p_j(t,y)) ,\eta(t)\Big)dt+\psi^{\sharp}(y).
\]
Recalling (\ref{PMP6}) and (\ref{u-eta}), the derivative of the cost w.r.t.~the terminal point of the trajectory is computed by
\begin{equation*}
\begin{split}
D\Phi_j(x)D_{y}x_j(0,y)&=~\int_{0}^{T}L_x\big(x_j,u^{\eta},\eta\big)  D_{z}x_j+ L_u\big(x_j,u^{\eta},\eta\big) {d\over dy}u^{\eta}dt+D\psi^{\sharp}(y)\\
&=~\int_{0}^{T}-{d\over dt}\left[p_j(t,y) D_{z}x_j(t,y)\right]dt+D\psi^{\sharp}(y)~
=~p_j(0,y)D_{z}x_y(0,y).
\end{split}
\end{equation*}
This implies\[
D\Phi_j(x)~=~D\Phi_j(x_j(0,y))~=~p_j(0,y)\,.
\]
By part (ii) of Lemma~\ref{l:42}, the a priori bound on (\ref{td2}) on the third derivative of 
$\psi^\sharp$ yields an a priori bound on the third derivative of the value function
$D^3 \Phi_j(x)$, for any $x= x_j(0,y)$ with $y\in \Gamma_j$.  Say,
\bel{tdg} \bigl| D^3 \Phi_j(x)\bigr|~\leq~M_j'\,.\eeq

%

We now apply part (i) Lemma~\ref{l:42}. This implies that, for any initial data (\ref{id3}), 
with $\|D^2 \vp\|_{\C^2}\leq  M'_j$, the solution to (\ref{ODE3}) satisfies a bound of the form
\bel{bp4}
\bigl|D^2_x p(T,x)\bigr|~\leq~M_j''.\eeq
The constant $M_k$ is now chosen so that
\bel{Mkdef} M_k~>~\max\{ M_1'', \ldots, M_{k-1}''\}.\eeq
We observe the above construction achieves the following:

Consider two families of trajectories satisfying the PMP,  starting in a neighborhood of the same point $x_0$,
and ending in different cubes, say $\Gamma'_j$ and $\Gamma'_k$, with $j<k$.
By the choice of $M_k$ at (\ref{Mkdef}) and the bounds (\ref{tdg}), at all initial points $y$
such that  $x_k(T, y)\in \Gamma_k'$, we have 
$$ \bigl| D^3 \Phi_j(x)\bigr|~\leq~M_j',\qquad\qquad  \bigl| D^3 \Phi_k(x)\bigr|~>~M_j'\,.$$
Indeed, if the second inequality did not hold, then we would have the bound (\ref{bp4}),
contrary to the construction of $\psi^\sharp$.

  Thus, the third derivatives $D^3\Phi_j(x)$ and $D^3\Phi_k(x)$ are strictly different 
  in a neighborhood of $x_0$. 
  \v
{\bf 6.} Based on the previous analysis, we give a bound on the Lebesgue measure
of the set of initial points $x_0$ from which two distinct optimal trajectories
initiate, ending in different cubes $\Gamma_j, \Gamma_k$.
This set contains:
\begi
\item Points $x_0 = x^\eta(0,y)$ with  $y\in B(0, R_0)$ such that the determinant
of the Jacobian matrix $D_y x^{\eta}(0,y)$ is small:
$$\bigl| \det \big(D_y x^\eta(0,y)\big)\bigr|~ \leq ~\delta_1\,.$$
 The Lebesgue measure of this set 
is  $< \delta_1 \cdot\meas\bigl(B(0, R_0)\bigr)$.    Choosing $\delta_1$ small enough, 
since the probability measure $\mu_0$ is absolutely continuous, we achieve (\ref{jsm}).

\item Points $x_0\in B(0, R_0)$ such that $x_0= x^\eta(0,y)$ for some
$y\in \Gamma_k\setminus \Gamma'_k$.   
By (\ref{gpri}) it follows 
$$\meas\left( \bigcup_{k=1}^\nu (\Gamma_k\setminus\Gamma'_k)\right)~<~\ve_1\,.$$
Again, since $\mu_0$ is absolutely continuous, by choosing $\ve_1>0$ sufficiently small,
we achieve
\bel{musm}
\mu_0\left( \Big\{ x^\eta(0,y)\,;~~y\in \bigcup_k(\Gamma_k\setminus\Gamma'_k)\Big\}\right)~<~{\ve_0\over 2}\,.\eeq
Toward (\ref{musm}), it is important to observe that the determinant of the Jacobian matrix  $D_y x^\eta(0,y)$ 
satisfies a uniform bound, depending on
the second derivatives $D^2\psi^\sharp$.   By (\ref{c2s}) these remain bounded, even when the third derivatives
are changed.

\item  The remaining set $S$ of all points $x_0\in B(0, R_0)$ which lie outside the previous two sets.
We claim that  $S$ has measure zero.  Indeed, if $x_0\in S$ is the initial point 
for two trajectories satisfying the PMP and terminating inside two distinct sets $\Gamma_j', \Gamma_k'$, then the corresponding value functions $\Phi_j, \Phi_k$ has distinct third derivative
 at $x_0$.  Therefore, $x_0$ cannot be a Lebesgue point of the coincidence set
 $\{ x\,;~\Phi_j(x)=\Phi_k(x)\}$.  Since the set has no Lebesgue points, it has measure zero.
 By the absolute continuity of $\mu_0$, we obtain 
 \bel{m0s}\mu_0(S)=0.\eeq
\endi
Combining the three bounds (\ref{jsm}), (\ref{musm}), and (\ref{m0s}), this achieves the proof.
\endproof

\section{Examples of structurally stable solutions}
\label{s:2}
\setcounter{equation}{0}
In this section we give  some examples of first order mean field games with one or more solutions,
and discuss their stability.

To motivate the examples concerning differential games, we first consider
 two maps of the unit disc $B_1\subset\R^2$ onto itself,
in polar coordinates $(r,\theta)$.
\bel{ex11}
\phi_1(r,\theta)~=~\left({2r\over 1+r^2}\,,~ \theta+\theta_0\right),\qquad\qquad 
\phi_2(r,\theta)~=~\left({r \over 1+r^2}\,,~ \theta\right)
,\eeq
where the rotation angle satisfies $0<\theta_0<2\pi$. 
Notice that the origin is  the unique fixed point of both $\phi_1$ and $\phi_2$. However, this fixed point is
asymptotically stable for the map $\phi_2$, but unstable for $\phi_1$. Indeed, for every $\bar r\geq 0$,  the sequence 
of radii
$$
r_{n+1}~=~{r_n\over 1+r_n^2},\qquad r_0~=~\bar r,
$$
is decreasing and converges to $0$.
On the other hand, for $0<\bar r<1$,  the sequence 
$$
r_{n+1}~=~{2r_n\over 1+r_n^2},\qquad r_0~=~\bar r,
$$
is increasing and converges to $1$.

\subsection{Games with a unique solution, stable or unstable.}
In the following examples of mean field games, 
as probability space labeling the various players we simply take $\Omega = [0,1]$. 
Motivated by (\ref{ex11}),  we begin by constructing  mean field games with a unique solution, 
which is unstable in the first example, and stable in the second.

\begin{example}{\rm  
Consider a game where each player $\xi\in [0,1]$ minimizes the same cost
\bel{dg7}
J(u)~=~\int_0^T |u(t)|^2\, dt + \bigl|x(T) - \psi(b(T))\bigr|^2,\eeq
subject to the trivial dynamics
\bel{ga2}\dot x(t) ~=~u(t),\eeq
with initial data
\bel{ga3} x(\xi,0)~=~\bar x(\xi)~=~0\qquad\forall \xi\in [0,1].\eeq
Here $u(t), x(t)\in\R^2$ while, as
in (\ref{bar}), $b(T)\in \R^2$ denotes the barycenter of the terminal positions of all players. 
Two cases will be considered.

\indent{\bf 1 - An unstable game.}  Let the terminal cost  be
\[
\psi(x)~=~{1+ T\over  T}\cdot \phi_1(x),
\]
where $\phi_1$ is the first map defined  at (\ref{ex11}), using polar coordinates. 
In this case,  $x(t, \xi )\equiv 0$ for all $(t,\xi)\in [0,T]\times [0,1]$ provides the unique solution to the mean field game. Indeed, given a barycenter $b(T)$, the PMP
\bel{PMP-EX11}
\left\{
\bega{rl} \dot x&=-\ds{p\over 2}\,,\\[3mm]
\dot p&=~ 0,\enda
\right.
\qquad\mathrm{with}\qquad
\begin{cases}
x(0)&=0\,,\\[3mm]
p(T)&=2(x(T)-\psi(b(T))),
\end{cases}
\eeq
has a unique solution 
\bel{uni}
x(t)~=~{ t\over 1+ T}\cdot \psi(b(T))\qquad t\in [0,T].
\eeq
All the optimal trajectories $x(\cdot,\xi)$ of the mean field game are the same.  In particular, if $x^*(t,\xi)$ is a solution to the game then
\[
b^*(T)~=~\int_{0}^1x^*(T,\xi)d\xi~=~x^*(T,\xi)~=~{ T\over 1+ T}\cdot \psi(b^*(T))~=~\phi_1(b^*(T)).
\]
Notice that $\phi_1$ has a unique fixed point, i.e.~the origin, we have  $b^*(T)=0$  and (\ref{uni}) yields  $x^*(\cdot, \xi )\equiv 0$ for all $\xi\in [0,1]$. On the other hand, for any sequence $b^{(k)}$ such that $b^{(k+1)}=\Phi(b^{(k)})$, one has that
$$b^{(k+1)}(T)~=~{\kappa T\over 1+\kappa T}\cdot \psi\big(b^{(k)}(T)\big)~=~\phi\big(b^{(k)}(T)\big).$$
Since $0$ is an unstable equilibrium of $\phi$,  the   zero solution of game is unstable.
\v

\indent{\bf 2 - A stable game.}  Similarly, if the terminal cost  $\psi$ is given by 
\[
\psi(x)~=~{1+ T\over  T}\cdot \phi_2(x)
\]
with $\phi_2$ being the first map in (\ref{ex11}) then  $x^*(\xi,\cdot)\equiv 0$ for all $\xi\in [0,1]$ is again the unique solution of the MFG. Moreover, since $0$ is asymptotically stable for the map $\phi_2$, the solution $x^*$ is stable. 
}
\end{example}

\subsection{Games with multiple solutions.}
Next, 
we give an example of a mean field game which admits both stable and unstable (but structurally stable) 
solutions. 


\begin{example}\label{e:21}{\rm   Here all  controls and trajectories are scalar functions. 
 The objective of every player is
\bel{ga1}
\hbox{minimize:} \quad \int_0^T \left[|u(t)|^2+ {1 \over 1+ x^2(t)} + \kappa\cdot
\bigl| x(t) - b(t)\bigr|^2\right]\, dt,\eeq
subject to
\bel{ga22} \dot x~=~u,\qquad\qquad x(0,\xi)~=~0\qquad\forall \xi\in [0,1].\eeq
Here  $b$ denotes the barycenter of the distribution of players as in (\ref{bar}).
\begin{proposition} For the MFG (\ref{ga1})--(\ref{ga22}), the following holds.
\begi
\item[(i)]
For all  $\kappa>1$ and  $T>2$, the mean field game has at least three  solutions.  These have the form
\bel{x012} x_i(t,\xi)~=~y_i(t),\qquad\qquad i=0,1,2,\eeq 
with $y_0(t)=0$, while $y_1$ is monotone increasing, and $y_2(t)= - y_1(t)$ for $t\in [0,T]$.
\item[(ii)] The zero solution is unstable.     However, assuming that $T\not={ (2n-1)\pi\over 2}$ for every 
$n\geq 1$, this solution is structurally stable.
\item[(iii)] Both  solutions $x_1, x_2$ are stable, and structurally stable.
\endi
\end{proposition}
{\bf Proof.}  {\bf 1.} Given a function $b(\cdot)$, the reduced Hamiltonian of (\ref{ga1})-(\ref{ga22})  is  computed by 
\[
\Hat H^{b}(x,p,t)~=~{1 \over 1+ x^2} + \kappa\cdot
\bigl| x - b(t)\bigr|^2-{p^2\over 2}.
\]
Assume that $\kappa >1$. For every $p\in\R,t\geq 0$, we have 
\[
\Hat H_{xx}^{b}(x,p,t)~=~2\cdot\left(\kappa-{1\over (1+x^2)^2}+{4x^2\over (1+x^2)^3}\right)~>~0\qquad\forall x\in\R,
\]
hence the map $x\mapsto \Hat H^{b}(x,p,t)$ is strictly convex. Thus, by the same argument in  Step {\bf 2} of the proof of Theorem \ref{t:31}, the optimal control problem (\ref{ga1})--(\ref{ga3}) has a unique optimal solution and  all the optimal trajectories 
 $x( \cdot,\xi)$ of the mean field game coincide. 
  As a consequence,  $x(\cdot,\xi)= b(\cdot)$ is an optimal solution to   the optimization problem
\bel{opt1}
\hbox{minimize:} \quad \int_0^T \left[{|\dot x|^2\over 2}+ {1 \over 2\bigl(1+ x^2(t)\bigr)}\right]\, dt,\quad
\hbox{subject to}\quad  x(0)~=~0.\eeq
Here, we can think of 
  $$K(\dot x)~=~{\dot x^2\over 2}\,,\qquad\qquad V(x)~=~-{1 \over 2(1+ x^2)}$$
respectively as kinetic and potential energy.  The solution is a motion governed by  the Euler-Lagrange equations
%
%
%
\bel{eq-tb1}
\ddot{y}(t)\,=\,- V_y(y)\,=\,-{ y(t)\over \big[1+ y^2(t)\big]^2},\qquad \qquad y(0)\,=\,0, \qquad \dot{y}(T)\,=\,0.
\eeq
It is clear that $y\equiv 0$ is a solution of (\ref{eq-tb1}) and this provides the first solution  of the mean field game
\bel{zeros} x(t,\xi)~=~0\qquad \forall \xi\in [0,1],\quad t\in [0,T].\eeq
To complete this step, we claim that  (\ref{eq-tb1})  admits at least  two additional solutions $y_1(\cdot)$, $y_2(\cdot)$,
  such that $y_1$ is  strictly increasing in $[0,T]$, and $y_2(t)= -y_1(t)$ . The mean field game has
 two more solutions $x_1, x_2$, as in (\ref{x012}).

 Observe that solutions to the Euler-Lagrange equations conserve the total energy
\bel{en} E(y, \dot y)~=~K(\dot y) + V(y)~=~{\dot y^2\over 2} - {1\over 2(1+y^2)}\,.\eeq
Level sets where $E$ is constant are plotted in Fig.~\ref{f:xepot}. Solutions to the boundary value problem
 (\ref{eq-tb1}) correspond to trajectories that start at time $t=0$ on the vertical axis where $y=0$, and end 
 at time $t=T$ on the horizontal axis where $\dot y=0$.

\begin{figure}[ht]
\centerline{\hbox{\includegraphics[width=10cm]{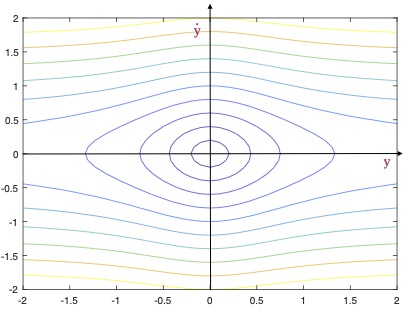}}}
\caption{\small  The level sets where the energy $E(y,\dot y)$ at (\ref{en}) is constant.}
\label{f:xepot}
\end{figure}

We thus seek an increasing solution of
\bel{yM}
\dot{y}(t)~=~ \sqrt{{1\over 1+ y^2(t)}-{1\over 1+ M^2}},\qquad \qquad y(0)~=~0,\eeq
for some constant $M$ such that $M=y(T)$. Calling $y=y(t,c)$ the solution to (\ref{yM}) with $M=c$, we have
\[
{\sqrt{c}\over 1+c^2}\cdot \sqrt{c-y(t,c))}~\leq~\dot{y}(t,c)~\leq~\sqrt{{2c\over 1+c^2}}\cdot \sqrt{c-y(t,c))}.
\]
By a comparison argument, we obtain for all $0\leq t\leq \sqrt{2(1+c^2)}$ that 
\bel{cm-y(t,M)}
c-c\cdot \left(1-{ t\over 2(1+ c^2)}\right)^2~\leq~y(t,c)~\leq~c-c\cdot \left(1-{ t\over \sqrt{2(1+ c^2)}}\right)^2.
\eeq
In particular, assume that $T>\sqrt{2}$. For every $c\geq \sqrt{{T^2-2\over 2}}$, the solution $y(\cdot,c)$ is defined on $[0,T]$ and satisfies  
\[
c\cdot \left(1-\left(1-{1\over T}\right)^2\right)~\leq~y(T,c)~\leq~ c.
\]
Calling
$M\doteq \inf\left\{c\geq \sqrt{{T^2-2\over 2}}: y(T,c)\leq c \right\}>0$, we claim that $y(T,M)=M$. 
Indeed, assume that $M-y(T,M)=\delta_0>0$. Then, by  (\ref{cm-y(t,M)}), one has   
\[
M-\delta_0~=~y(T,M)~\leq~M\cdot\left(1-\left(1-{ T\over \sqrt{2(1+ M^2)}}\right)^2\right).
\]
Hence, $M -\ds \sqrt{T^2-2\over 2}=\ve_0>0$ and the map $t\mapsto y(t,M-\ve)$ is defined on $[0,T]$ for all $0< \ve<\ve_0$. Moreover, by the monotone increasing property of $c\mapsto y(T,c)$, we have 
\[
y(T,M-\ve)~\leq~y(T,M)~=~M-\delta_0~\leq~M-\ve
\]
for all $0<\ve<\min\{\ve_0,\delta_0\}$. This yields a contradiction.

In the next steps we will show that  all three solutions are essential, the zero solution is unstable,  and the  two non-zero solutions are stable.
\medskip

{\bf 2.} 
We begin by showing that the null solution $x(t,\xi)\equiv 0$  is unstable but essential. In the present case,  the map $b\mapsto \Tilde b=\Phi(b)$ at (\ref{Phi})-(\ref{mot})  takes the form 
\[
\Phi(b)(t)~=~x_b(t)\qquad\forall t\in [0,T],
\]
where $(x_b,p_b)$ denotes the unique solution of the PMP
\bel{PMP-EX1}
\left\{
\bega{rl} \dot x&=~u(x,p)~=~-\ds{p\over 2},\\[3mm]
\dot p&=~ \ds{2 x\over (1+ x^2)^2}- 2\kappa\cdot (x-b),\enda
\right.
\qquad
\left\{
\bega{rl} x(0)&=~0,\\[3mm]
p(T)&=~ 0,\enda
\right.
\eeq
where the optimal control is 
\[ u(x,p)~=~\argmin_{\omega\in\R} \Big\{ \omega^2 + p\,\omega\Big\}~=~-{p\over 2}.
\]

Linearizing the system  (\ref{PMP-EX1}) at $b\equiv 0$ we obtain an expression for the 
differential $D\Phi(0)$, namely
$$D\Phi(0) \, b~=~\Hat b,$$
where $\Hat b(t) = X(t)$ is the function obtained by solving the linear system
\bel{Li-PMP-EX1}
\begin{bmatrix}
\dot X(t)\\[4mm]
\dot P (t)
\end{bmatrix}~=~\begin{bmatrix}0& -1/2\\[4mm]
2-2\kappa& 0 
\end{bmatrix}
\begin{bmatrix}
X(t)\\[4mm]
P(t)
\end{bmatrix}+2\kappa b(t) \begin{bmatrix}0\\[4mm] 1
\end{bmatrix},
\qquad X(0)~=~P(T)~=~0.
\eeq
%
%
%
%
Eliminating the variable $P= -2\dot X$, one is led to the second order ODE
$$ -2 \ddot Y ~=~ (2-2\kappa) Y + 2\kappa  b.$$
To determine eigenvalues $\lambda$ and eigenfunctions $Y$, we need to solve
$$ -2 \ddot Y = (2-2\kappa) Y + {2\kappa\over \lambda} \,Y, \qquad\qquad Y(0)\,=\, \dot Y(T)\,=\,0.
$$
\bel{EE}  \ddot Y + \left(1-\kappa+{\kappa\over\lambda}\right) Y ~=~0, \qquad\qquad Y(0)\,=\, \dot Y(T)\,=\,0.
\eeq
The eigenvalues and eigenfunctions of $D\Phi(0)$ are thus found to be
\bel{EGV}
\lambda_n~=~
{\kappa\over \kappa+{(2n-1)^2\pi^2\over 4T^2 } -1}\,,\qquad\qquad Y_n(t)~=~\sin \left({(2n-1)\pi\over 2T}t\right), \qquad n=1,2, \dots.
\eeq
In particular, if $T>\ds {\pi\over 2}$ and $\kappa>1$, computing  the first eigenvalue of $D\Phi(0)$ one finds
 $\lambda_1 > 1$. This implies that the null solution $x(t,\xi)\equiv 0$ is unstable.

On the other hand,  we observe that, by (\ref{EGV}), if
\bel{Tne}
T\,\not=\, {(2n-1)\pi\over 2}\qquad\qquad\hbox{for every} ~ n\geq 1,\eeq then  $1$ is not an eigenvalue of $D\Phi(0)$. In this case, using  the same argument as in Step {\bf 4} of the proof of Theorem \ref{t:31}, 
we conclude  that $y_1$ is essential.
\v

{\bf 3.} We now prove that $y_1$ is stable.  Given any ${\bf b}\in \mathcal{C}([0,T])$, we first compute $D\Phi(y_1)(\bar{b})$. As in step 2, for every $\ve\in\R$, let  $(x^{\ve}(t), p^{\ve}(t))$ be the solution of $(\ref{PMP-EX1})$ with $b=y_1+\ve\bar{b}$. By the linearization, it holds 
\[
\begin{bmatrix}
x^{\ve}(t)\\[4mm]
p^{\ve}(t)
\end{bmatrix}
~=~\begin{bmatrix}
y_1(t)\\[4mm]
p_1(t)
\end{bmatrix}+\ve \,\begin{bmatrix}
{\bf x}_{\bar{b}}(t)\\[4mm]
{\bf p}_{\bar{b}}(t)
\end{bmatrix}+o(\ve).
\]
Here $\begin{bmatrix}
{\bf x}_{\bar{b}}(t)\\
{\bf p}_{\bar{b}}(t)
\end{bmatrix}$ is the solution to  the equation obtained linearizing (\ref{PMP-EX1})  around $y_1$, namely
\bel{LN-PMP-y1}
\begin{cases}
\dot{\bf x}(t)~=~-\ds{{\bf p}(t)\over 2}\,,\\[3mm]
\dot{\bf p}(t)~=~\ds 2\left({1-3y_1^2\over (1+y_1^2)^3}-\kappa \right){\bf x}+2{\kappa}\bar{b}\,,
\end{cases}
\qquad \begin{cases}
{\bf x}(0)~=~0,\\[3mm]
{\bf p}(T)~=~0.
\end{cases}
\eeq
Let  the pair $(\gamma,\bar{b})$ denote an eigenvalue and an eigenfunction of $D\Phi(y_1)$. As in Step {\bf 2}, we have 
\[
D\Phi(y_1)(\bar{b})~=~{\bf x}_{\bar{b}}~=~\gamma \bar{b},
\]
and ${\bf x}_{\bar{b}}$ solves the two point boundary problem 
\bel{2-p-b}
\ddot{y}(t)~=~\left[\kappa\cdot\left(1-{1\over \gamma}\right)+{3y_1^2-1\over (1+y^2_{1})^3 }\right]\cdot y(t),\qquad y(0)~=~\dot{y}(T)~=~0.
\eeq
To verify the stability of $y_1$, we will show that  all  eigenvalues of $D\Phi(y_1)$  are contained within the open
interval $]0,1[$. Assume by a contradiction that $D\Phi(y_1)$ has an eigenvalue $\gamma\in \R\setminus\,]0,1[\,$, 
so that the equation (\ref{2-p-b}) has a nonzero solution $y_2$. Recalling that $t\mapsto y_1(t)\in [0,+\infty[$ is strictly increasing with $y_1(0)=0$, we define 
\[
t_1~\doteq~\min\left\{t\in [0,T]: \kappa\cdot \left(1-{1\over \gamma}\right)+{3y_1^2(t)-1\over (1+y^2_{1})^3(t) }\geq 0\right\}.
\]
For every $\tau\in [t_1,T]$,  from (\ref{2-p-b}) it follows
\[
-y_2(\tau)\dot{y}_2(\tau)~=~\int_{\tau}^{T}\dot{y}_2^2(t)dt+\int_{\tau}^{T}\kappa\cdot\left(1-{1\over \gamma}\right)+{3y_1^2-1\over (1+y^2_{1})^3 } y_2^2(t)dt~>~0.
\]
Therefore, both $y$ and $\dot{y}$ do not change sign in $[t_1,T]$. Without loss of generality, we can assume that $y$ is positive in $[t_1,T]$. Set 
\[
t_2~\doteq~\max\bigl\{t\in [0,t_1]:y_2(t)~=~0\bigr\}.
\]
We then have 
\[
y_2'(t_2)~>~0,\qquad y_2(t_2)~=~0,\qquad\mathrm{and}\qquad y_2(t)~\geq~0\qquad\forall t\in [t_2,T].
\]
On the other hand, since $y_1$ is an increasing solution of (\ref{eq-tb1}), the function $z_1\doteq \dot{y}_1$ solves the equation
\[
\ddot{z}(t)~=~{3y_1^2-1\over (1+y^2_{1})^3 }\cdot z(t),~~~\qquad \dot{z}(0)~=~z(T)~=~0.
\]
Thus, for all $t\in [t_2,T]$, one has 
\[
 [\dot{y}_2(t)z_1(t)]'~=~ [\dot{z}_1(t)y_2(t)]'+\left(1-{1\over \gamma}\right) y_2(t)z_2(t)~\geq~[\dot{z}_1(t)y_2(t)]',
\]
and this yields
\bel{ode-y2}
\dot{y}_2(t)z_1(t)-\dot{y}_2(t_2)z_1(t_2)~\geq~\dot{z}_1(t)y_2(t).
\eeq
Equivalently, 
\[
{d\over dt} \left({y_2(t)\over z_1(t)}\right)~\geq~\dot{y}_2(t_2)z_1(t_2)\cdot {1\over z_1^2(t)}\qquad\forall t\in [t_2,T].
\]
This implies
\[
y_2(t)~\geq~\dot{y}_2(t_2)z_1(t_2)\cdot \int_{t_2}^{t}{1\over z_1^2(t)}dt.
\]
Therefore, by (\ref{ode-y2}) one has
\bel{dot-y2}
\dot{y}_2(t)~\geq~\dot{y}_2(t_2)z_1(t_2) \cdot\left[{1\over z_1(t)}+\dot{z}_1(t)\cdot \int_{0}^{t}{1\over z^2_1(s)}ds\right].
\eeq
To obtain a contradiction, we will show that 
\bel{check}
0~=~\dot{y}_2(T)~=~\dot{y}_2(t_2)z_1(t_2) \cdot\lim_{t\to T-}\left[{1\over z_1(t)}+\dot{z}_1(t)\cdot \int_{0}^{t}{1\over z^2_1(s)}ds\right]~>~0.
\eeq
Assume that $y_1(0)=v_0$ and $\beta=y_1(T)$.  We then have 
\[
z_1(t)~=~\dot{y}_1(t)~=~\left(v_0^2-{y^2_1(t)\over 1+y^2_1(t)}\right)^{1/2},\qquad 
v_0^2~=~{\beta^2\over 1+\beta^2}\,,
\]
and 
\[
\dot{z}_1(t)~=~\ddot{y}_1(t)~=~-{y_1(t)\over (1+y^2_1(t))^2}.
\]
By a change of variable, (\ref{check}) is equivalent to 
\[
I~\doteq~\lim_{y\to \beta-} \left[{(1+\beta^2)^{1/2}(1+y^2)^{1/2}\over (\beta+y)^{1/2}(\beta-y)^{1/2}}-{y\over (1+y^2)^2}\cdot \int_{0}^{y}{(1+\beta^2)^{3/2}(1+z^2)^{3/2}\over (\beta+z)^{3/2}(\beta-z)^{3/2}}dz\right]~>~0.
\]
Notice that for $\beta>0$ sufficiently large, we have 
\[
\lim_{y\to \beta} \left[{(1+\beta^2)^{1/2}(1+y^2)^{1/2}\over (\beta+y)^{1/2}(\beta-y)^{1/2}}-{1+\beta^2\over \sqrt{2\beta}(\beta-y)^{1/2}}\right]~=~0\,,
\]
and 
\begin{eqnarray*}
{y\over (1+y^2)^2}\cdot\int_{0}^{y}{(1+\beta^2)^{3/2}(1+z^2)^{3/2}\over (\beta+z)^{3/2}(\beta-z)^{3/2}}dz&\leq&{y(1+\beta^2)^{3/2}\over (1+y^2)^{1/2}(\beta+y)^{3/2}}\int_{0}^{y} (\beta-z)^{3/2}dz\\
&=&{2y(1+\beta^2)^{3/2}\over (1+y^2)^{1/2}(\beta+y)^{3/2}}\cdot\left[{1\over (\beta-y)^{1/2}}-{1\over \beta^{1/2}}\right].
\end{eqnarray*}
In particular, this implies
\[
I~\geq~{1+\beta^2\over \beta\sqrt{2}}+(1+\beta^2)\cdot \lim_{y\to \beta-}\left({1\over \sqrt{2\beta}}-{2y(1+\beta^2)^{1/2}\over (1+y^2)^{1/2}(\beta+y)^{3/2}}\right)\cdot {1\over \beta-y}~=~{1+\beta^2\over \beta\sqrt{2}}~>0~.
\]
This shows that all  eigenvalues of $D\Phi(y_1)$  are contained in the open interval $]0,1[\,$, and $y_1$ is a stable solution of the MFG. By symmetry, 
$x_2(t,\xi)~=~y_2(t)\doteq -y_1(t)$ for $t\in [0,T]$ and all  $\xi\in [0,1]$, is also a stable solution of the MFG.
\endproof}
\end{example}

\subsection{Examples of games with no solutions.}

\begin{example}\label{e:22}{\rm 
Consider the mean field game on the time interval $t\in [0,T]$,
where player $\xi\in  \Omega = [0,1]$ has dynamics
\bel{ex}
\dot{x}~=~u~\in~[-1,1],\qquad x(0,\xi)~=~0. 
\eeq
%
The goal of player  $\xi$ is to optimize his terminal position relative to the distribution of the other players, namely
\bel{ex2}
\hbox{maximize:} \qquad  \bigl|x(T,\xi)- b(\xi)\bigr|^2,\eeq
where
\bel{ex3} 
b(\xi)~=~\int_0^1 e^{-|\zeta-\xi|^2}\cdot x(T,\zeta)\, d\zeta.\eeq

We claim that this game has no strong solution. Indeed, if $b(\xi)\equiv 0$,
then every player has two equally good strategies: 
\bel{ex4}u(t)\equiv 1, \qquad x_2(t) = t\qquad\hbox{or}\qquad u(t)\equiv -1,\qquad x_2(t) = -t.
\eeq
This cannot be a solution, because  $\xi\mapsto x(T,\xi)\in \{-T,T\}$ is a measurable map, and
the integral in (\ref{ex3}) cannot be identically zero.

On the other hand, if $b(\xi)$ is not identically zero, then
$$\int_0^1 b(\xi) \,x(T,\xi)\, d\xi~=~ \int_0^1 b(\xi) \cdot (-T \hbox{sign} \,b(\xi)\bigr)\, d\xi~=~
- T \int_0^1 \bigl| b(\xi)\bigr|\, d\xi ~<~0.$$
However, the definition of $b$ implies
$$\bega{rl}\ds
 \int_0^1 b(\xi) \,x(T,\xi)\, d\xi&\ds=~ \int_0^1 \left( \int_0^1 x_2(T,\zeta)e^{-|\zeta-\xi|^2} \, d\zeta\right) x(T,\xi)\, d\xi
 \\[4mm]
&\ds =~\ds\int_0^1\int_0^1 e^{-|\zeta-\xi|^2}\, x(T,\zeta)\, x(T,\xi)\, d\zeta d\xi~\geq~0,\enda
$$
%
reaching a contradiction.\footnote{Indeed, if the kernel can be written as the convolution
$\vp*\vp$, for some even function  $\vp(z)$, rapidly decreasing as $|z|\to \infty$, then
(replacing $z$ with $z-y$ as variable of integration and using the fact that $\vp(s)=\vp(-s)$)
$$\bega{l}\ds \dint (\vp*\vp)(x-y) f(x) f(y)\, dxdy~=~\int\!\!\int\!\!\int \vp(x-y-z)\vp(z) f(x)f(y)\, dz dx dy
\\[3mm]
\qquad=~\ds \int\!\!\int\!\!\int \vp(z-x)\vp(z-y) f(x)f(y)\, dz dx dy~=~\int (\vp*f)(z)\cdot (\vp*f)(z)\, dz~\geq~0.
\enda
$$}
Notice that here the unique mild solution is a measure, where each player uses the
two controls in (\ref{ex4}) with equal probability. 

}
\end{example}

\begin{example}\label{e:43}{\rm 
Consider the mean field game on the time interval $t\in [0,T]$,
where all players have the same dynamics and the same cost functional:
\bel{ce1}
\hbox{minimize:} \qquad\int_0^T u^2(t)\, dt + \psi\bigl( x(T)\bigr),\eeq
subject to
\bel{ce2}  \dot x~=~u - b^2,\qquad
\quad |u(t)|\leq 1,
\qquad\qquad
x(0,\xi)~=~0\qquad\forall \xi\in \Omega\,, \eeq
and with terminal constraint
\bel{ce3} \vp(x(T)\bigr)~\doteq~ \bigl(T- x(T)\bigr)\cdot x(T)~=~0.\eeq
Here 
\bel{ce4} b(t)~\doteq~\int_\Omega x(t,\xi)\, d\xi\eeq
denotes the barycenter of the distribution of players at time $t$,
while the terminal cost is a smooth function that satisfies
\bel{ce5}
\psi(x)~=~\left\{ \bega{cl} 0\quad & \hbox{if}\quad x=0,\cr
-2T\quad& \hbox{if}\quad x=T.\enda\right.\eeq
Notice that  the terminal constraint (\ref{ce3}) is equivalent to
\bel{ce6} x(T)~\in~\{0,T\}.\eeq

We claim that this mean field game has no solution. Namely,  
the ``best reply map" $\bfX\mapsto \Psi(\bfX)$ 
from $\L^1\Big( \Omega\,;~\C\bigl([0,T];\, \R^n\bigr)\Big)$ into itself 
does not have
any fixed point. To prove this, 
consider first the case where $\bfX = {\bf 0}\in \L^1\Big( \Omega\,;~\C\bigl([0,T];\, \R^n\bigr)\Big)$. That means:
\bel{zero} x(t,\xi)~=~0\qquad \hbox{for all ~~ $t\in [0,T]$~~and ~~$\mu$-a.e.~$\xi\in\Omega$.}
\eeq
In this case, $b(t)=0$ for all $t\in [0,T]$. Hence the optimal strategy for every player is
to choose
$u(t,\xi)=1$.   The corresponding trajectory $x(t,\xi)=t$ satisfies the terminal 
constraint  (\ref{ce3}) and achieves minimum cost 
$$J_{\min} ~=~\int_0^T 1\, dt +\psi(T) ~=~T-2T~=~-T.$$
On the other hand, if (\ref{zero}) fails, then
$b(t)$ is not identically zero and  the solution to 
(\ref{ce2}) cannot attain the value 
$x(T)= T$.   Hence the best strategy for every player is to take $u(t,\xi)\equiv 0$,
which yields the trajectory $x(t,\xi)=0$, with zero cost.

We have thus shown that
$${\bf 0}~ \notin~\Psi({\bf 0}), \qquad \hbox{while}\quad \Psi(\bfX)~=~\{{\bf 0}\}
\quad \forall ~\bfX\not= {\bf 0},$$
hence $\Psi$ cannot have a fixed point.

Notice that in this example the mean field game does not even admit mild solutions, in the randomized sense.}
\end{example}

We observe that in this example, the minimum cost does not depend continuously on the 
parameter $b(\cdot)$.    Namely, it jumps from $0$ down to  $-T$ as $b$ becomes the zero 
function. This is due to a lack of transversality in connection with the terminal constraint.

\section{Concluding remarks}
\label{s:9}
\setcounter{equation}{0}
In this paper we considered a class of first order mean field games,
characterized by 5-tuples $(f,L,\psi, \phi,\bar x)$ specifying the dynamics, cost functionals, 
averaging kernels, and initial distribution of players.

The main results show that, generically, for every given $\eta(\cdot)$  a.e.~player has a unique
optimal control $t\mapsto u^\eta(t,\xi)$.  As a consequence, 
 the ``best reply" map  $\eta\mapsto \Phi(\eta)$ at (\ref{mot}) is 
single valued, and the MFG has a strong solution.
Moreover, there are open sets of games with unique solutions, and open sets of games with multiple solutions.
These can be stable, or unstable, in the sense of Definition~\ref{d:22}.

It would be of interest to analyze whether similar results remain valid in a more general setting. Namely:
\begi
\item[(i)] Systems with fully nonlinear dynamics, i.e.~where the function $f(x,u,\eta)$
in (\ref{mdyn})  is not necessarily
affine w..r.t.~the control. 
\item[(ii)] Optimal control problems in the presence of terminal constraints, say
$$g_i\bigl(x(T,\xi)\bigr)~=~0,\qquad\qquad i=1,\ldots,N.$$
\endi

In all our previous examples, the mean field games had structurally stable solutions.
We thus  conclude the paper with a natural conjecture:

\begin{conjecture} 
For a generic 5-tuple $(f,L, \psi,\phi,\bar x)\in \X$, the MFG (\ref{etai})--(\ref{iitt}) has 
finitely many solutions, all of which are structurally stable.
\end{conjecture}

\v

\end{document}